\documentclass[12pt]{article}
\usepackage{amsmath,amsfonts,amsthm,amssymb}
\usepackage{epsfig}
\newcommand{\eps}{\varepsilon}

\newcommand{\mat}{\mathbf{M}}

\newtheorem {thm}     {Theorem}
\newtheorem {lem} {Lemma}
\newtheorem {rmk}  {Remark}
\newtheorem {cor} {Corollary}
\newtheorem {deff} {Definition}
\newtheorem {assump} {Assumption}
\newtheorem{prop}{Proposition}

\def\({\begin{eqnarray}}
\def\){\end{eqnarray}}
\def\[{\begin{eqnarray*}}
\def\]{\end{eqnarray*}}

\def\R{\mathbb{R}}

\def\d{{\mathrm{d}}}
\def\e{{\mathrm{e}}}
\def\O{{\mathcal{O}}}

\def\bydef{:=}

\title{The continuous limit of the Moran process and the diffusion
of mutant genes in infinite populations}

\author{Fabio A. C. C. Chalub\thanks{%
Departamento de Matem\'atica and Centro de Matem\'atica
e Aplica\c c\~oes, Universidade Nova de Lisboa, 
Quinta da Torre, 2829-516, Caparica, Portugal.
e-mail:{\tt chalub@cii.fc.ul.pt}},
Max O. Souza\thanks{%
Departamento de Matem\'atica Aplicada, Universidade Federal
Fluminense, R. M\'ario Santos Braga, s/n, 22240-920, Niter\'oi, RJ, Brasil.
e-mail:{\tt msouza@mat.uff.br}}}

\date{\today}

\begin{document}

\maketitle

\begin{abstract}
We  consider the  so called Moran process with frequency dependent
fitness given by a certain pay-off matrix. For finite populations, we
show that the final state must be homogeneous, and show how to compute
the fixation probabilities. Next, we consider the infinite population
limit, and discuss the appropriate scalings for the drift-diffusion
limit. In this case, a degenerated parabolic PDE is formally obtained
that, in the special case of frequency independent 
fitness, recovers
the celebrated Kimura equation in population genetics.
We then show that the corresponding initial value
problem is well posed and that the discrete model converges to the PDE
model as the population size goes to infinity. We also study some
game-theoretic aspects of the dynamics and characterize the best
strategies, in an appropriate sense.
\end{abstract}

\section{Introduction}

Since the beginning of modern evolutionary theory, the study
of the dynamics of a mutant gene in a population has 
attracted attention~\cite{Fisher1,Fisher2,Haldane,Wright2,Wright1}. 
It has been known for a long time that a mutant gene will
be, eventually, either fixed or lost.  The final 
result depends not only on natural selection but
also on chance~\cite{Kimura}.

The most natural attempt to describe mathematically the
evolution of a mutant gene uses a discrete model
for a finite population. The question of finding a consistent
model for the infinite population is then a natural one.
This is called in the physical literature the ``thermodynamical
limit'', and it is a classical subject on that field. See, for
example,~\cite{DombGreenVol1}.

When we consider the infinite population limit, it is natural
to have continuous
variables where we previously had discrete ones. For
example, if  $N$ denotes the size of the population,
the possible fractions of mutants are $0,1/N,2/N,\cdots,1$.
In the limit $N\to\infty$ this fills  the
entire real interval $[0,1]$. Time should
be rescaled accordingly, such that 
that the probability of fixation, for a given fraction $x$ over a
given time span $t$, does not depend significantly on the size
  of the population. In the infinite limit, we obtain a  partial
  differential equation (PDE).
This PDE is an approximation for
large $N$ of the discrete process, and as such
must present diffusion to the boundaries as a continuous 
representation of the fact that the mutant gene will be
eventually fixed or lost.

A different approach to the same problem is to consider continuous
models from the beginning. This leads to two distinct modeling
paradigms: the first one uses ordinary differential
equations (ODEs) to model the evolution of the fraction of mutants.
The most widely used equation in this context is the replicator
dynamics and its variations~\cite{SigmundHofbauer}. The second one,
which will be further developed in this work, uses PDEs to model
the evolutionary process. This approach goes back, at least, to the
seminal work by Kimura~\cite{Kimura}, where it was used to model
the diffusion of mutant genes. More explicitly, Kimura considered
the probability of fixation of a mutant gene that in a given time is
present in a certain fraction of the population. Here, we will deduce
the Kimura equation as a particular case of our work, where a PDE
will be obtained from the  more basic discrete process, in the
infinite population limit,  for the diffusion of mutant genes.

In this work, we consider a simple evolutionary process,
called the Moran process, introduced in~\cite{Moran} (used, e.g.,
for cancer dynamics~\cite{NowakFranziska2,NowakFranziska}, 
paleontology~\cite{Nee1},
phylogeny~\cite{Nee2}, genealogy~\cite{Donnelly1},
and epidemiology~\cite{Welch1}) for a
finite population and obtain a partial differential
equation as its thermodynamical limit. 

Our starting point
is evolutionary game theory. We consider
a finite population of fully connected 
interacting individuals through a certain
pay-off matrix. 
We start by proving that for any finite population (of size
$N$) of two-types, one of the types will be fixed after long 
enough time. The thermodynamical limit is then obtained
as a PDE that approximates the finite population dynamics for large $N$. 
We consider two different scalings for the time-step $\Delta t$, 
namely: $\Delta t=1/N$ (drift limit) and $\Delta t=1/N^2$
(drift-diffusion or simply diffusion limit). In the second
case it is also important 
 to introduce the so called weak selection
limit (pay-offs go to 1, when population goes to infinity).
We also show that the most interesting equations appear
in the drift-diffusion limit.

All the equations found in the limit are
degenerate, i.e., the diffusion coefficient vanishes on the
boundaries. 
The mathematical theory for such equations is not as well developed as
for the non-degenerate case. There are the classical books
\cite{CS76,DiBenedetto93}. In particular, \cite{CS76} proves existence
and uniqueness for the equation obtained by Kimura. 
For more recent works, see also~\cite{Dolb2,Dolb1}.

If we impose no diffusion in the PDE model, 
the solution can be decomposed in point dynamics, where each
  fraction evolves  through the replicator 
dynamics. The stationary states and long time behavior
of the replicator dynamics are, however, different to
the ones obtained as the thermodynamical limit of
the final states of discrete populations, showing that the
diffusion is essential to understand the discrete dynamics.

The PDE model allows the introduction of a relation
of dominance between
two different strategists that turns out to be, in its
dynamical features, identical to the flow 
of the replicator dynamics. We also show that the best
possible strategy in the finite, but large, population case
is given by the evolutionarily stable strategy 
(ESS) of the game~\cite{SigmundHofbauer,JMS}. 
This clarifies
the relation between a homogeneous population 
playing mixed strategies with given frequencies and 
a mixed population, with constant fractions, playing
pure strategies, i.e., the difference between evolutionarily
stable strategies and evolutionarily stable states. 

If the fitness for the individuals in the
population is frequency independent, the
resulting equation
is equivalent to a well-know equation of population
genetics, introduced by Kimura~\cite{Kimura}, describing
the probability of fixation of a mutant 
(with frequency independent fitness) in a population.
It turns out, that the equation derived in this work and the one introduced
by Kimura are a forward/backward pair of equations.

It is important to note that the perception that 
the Moran process (at least in the frequency independent
case) is related to diffusion process is not 
new~\cite{Iorio1,Iorio2}. The compatibility
between finite populations simulations and
the ESS, defined in the continuous case, are
also studied in~\cite{fogel1,fogel2,Orzack}.

The structure of this work is the following: 
In Section~\ref{replicator} we introduce the
(finite population) Moran process and study 
its properties. In particular, we prove that
the final state will be always homogeneous.
In Section~\ref{sec:thermo} we introduce the
drift-diffusion scaling and obtain a PDE
as the thermodynamical limit of the Moran
process. We also study its dynamic features from the
strategic point of view.
In Section~\ref{diffusionless}, we consider 
the no diffusion case and compare the
PDE obtained with the replicator dynamics.
In Section~\ref{moran} we particularize all
results to the frequency independent case
and in Section~\ref{drift} we study the
drift scaling. Finally, in Section~\ref{final}
we point new directions for this work, showing
how the tools developed here can be applied
to different dynamics.

\section{The frequency dependent discrete case}
\label{replicator}

\begin{figure}

\begin{center}
\epsfig{file=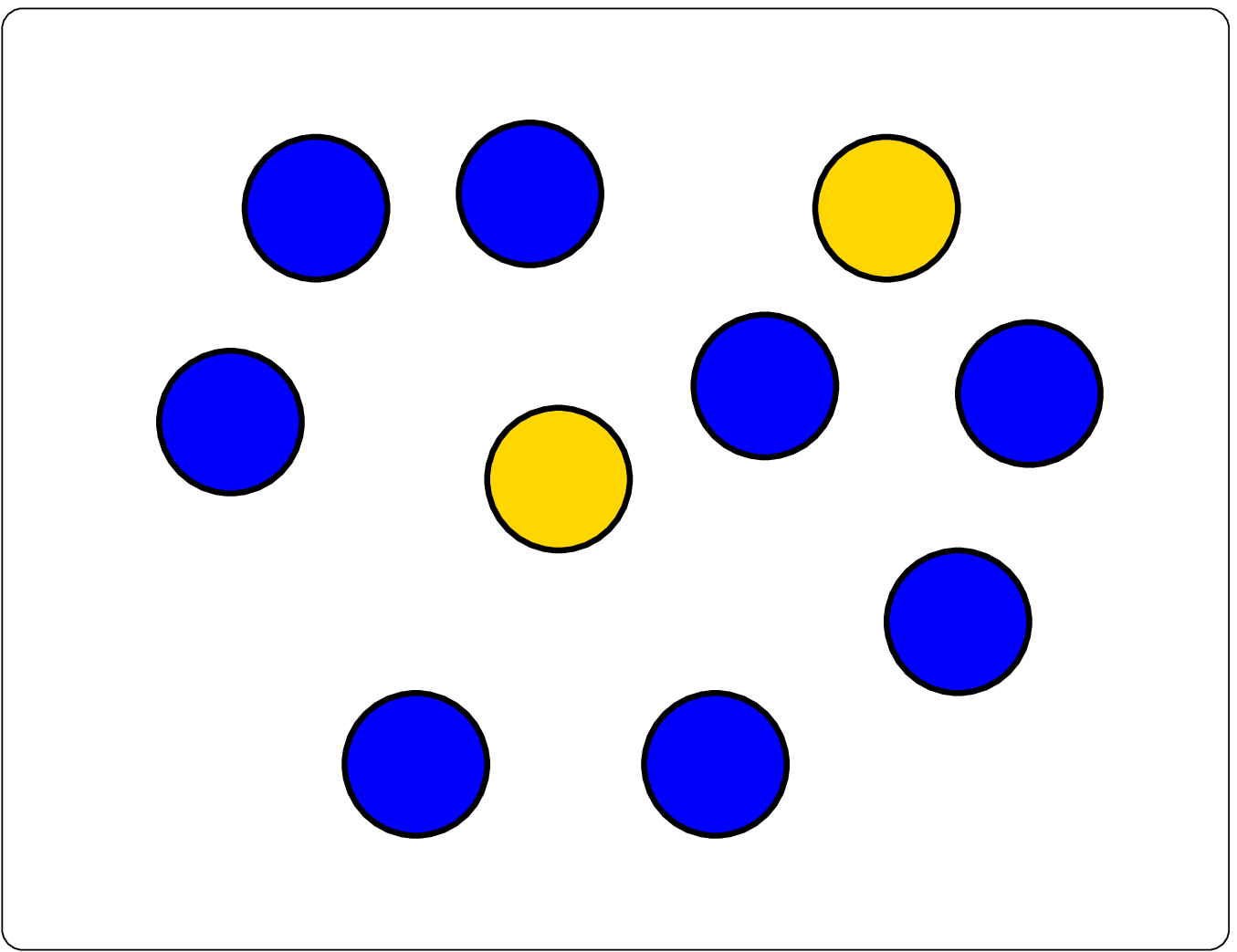,width=0.3\textwidth}
\hspace{.01\textwidth}
\epsfig{file=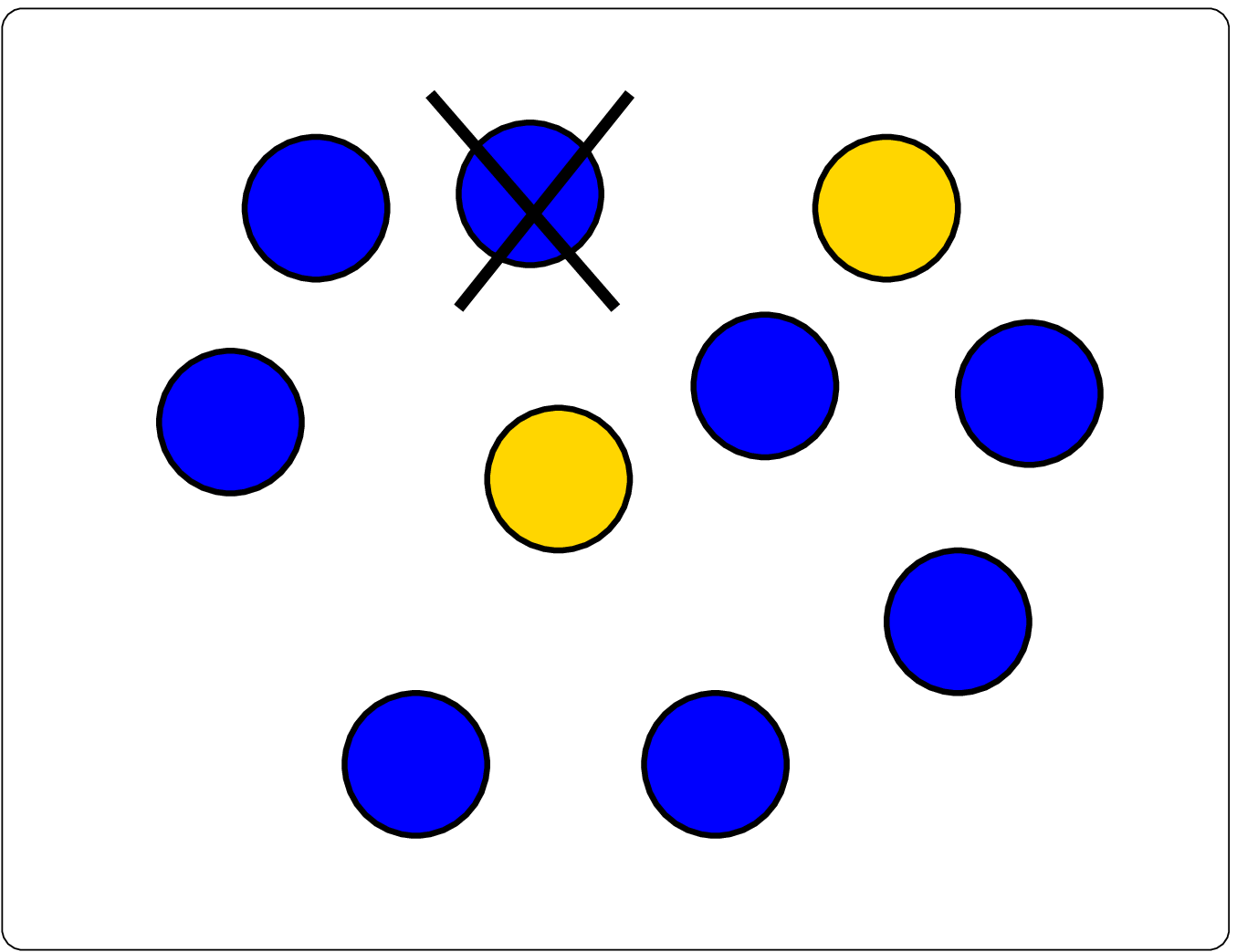,width=0.3\textwidth}
\hspace{.01\textwidth}
\epsfig{file=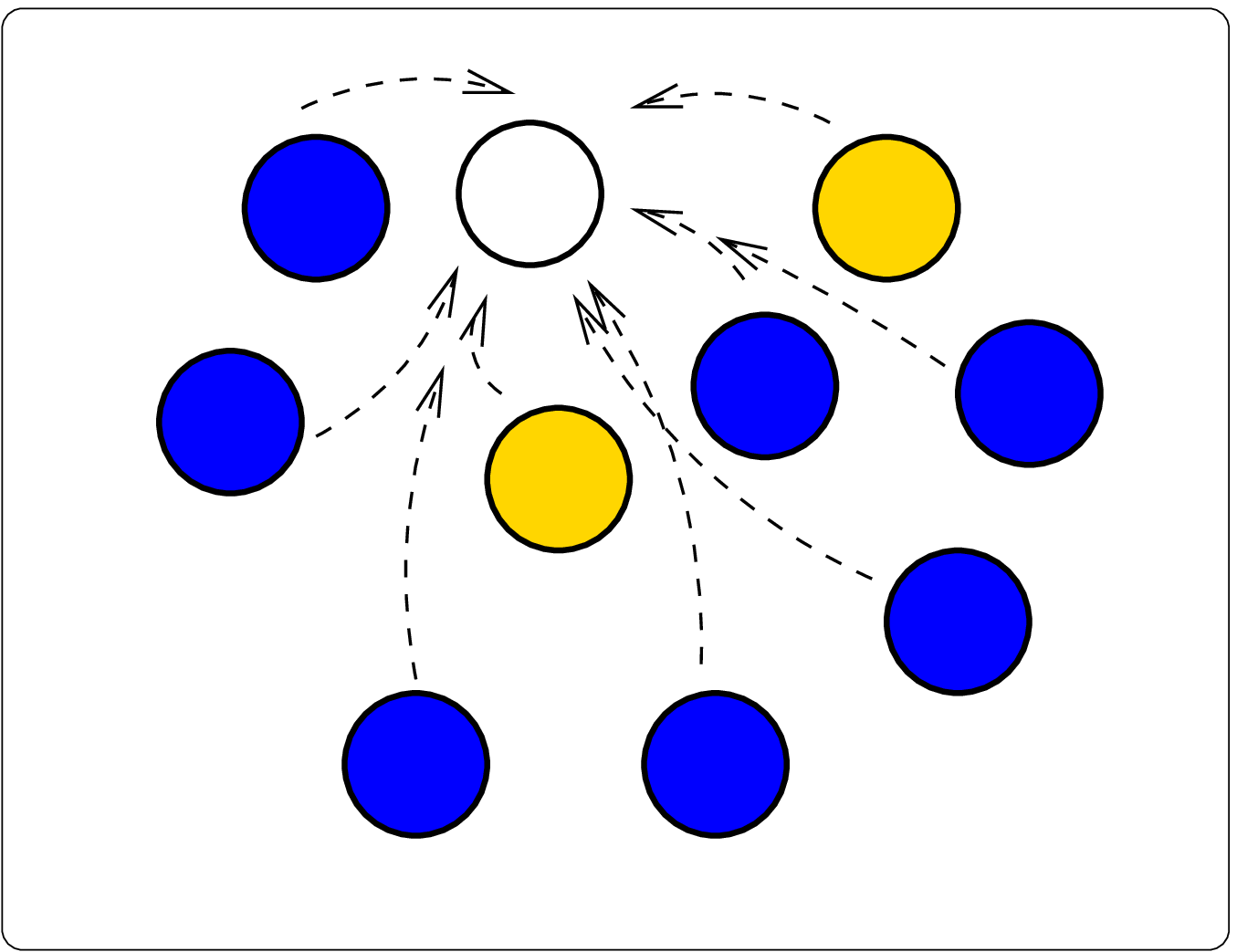,width=0.3\textwidth}
\vspace{.1cm}
\text{(a)\hspace{.29\textwidth} (b)\hspace{.29\textwidth} (c)}
\end{center}
\caption{The Moran process: from a two-types population (a) we chose 
one at random to kill (b) and a second to copy an paste in the place
left by the first, this time proportional to the fitness.
\label{Moran_fig}}
\end{figure}

We consider a fixed size population with two types of 
individuals: $\mathbb{A}$ and $\mathbb{B}$, say. At fixed
time steps, we choose one of the individuals to 
be eliminated at random and replace it by
a newborn which can be of either type. This newborn
is obtained as a copy of one of the remaining individuals with
probability proportional to its fitness. See Figure~\ref{Moran_fig}
for an illustration. This process is called the 
Moran process~\cite{Moran}.

Let $P(t,n,N)$ be the probability that there are $n$ type
  $\mathbb{A}$ individuals  at time $t$ in a population of fixed size $N$.
We define $c_+(n,N)$ ($c_0(n,N)$ and $c_-(n,N)$, respectively) 
as the probability
(independent of time) that the number of mutants
changes in time $t$ from $n$ to $n+1$ (to $n$ and to $n-1$ respectively)
in time $t+\Delta t$. We assume that these transition probabilities are
proportional to the fitness $\phi_A$ and $\phi_B$ of
types $\mathbb{A}$ and $\mathbb{B}$ respectively; thus we have:
\begin{eqnarray}
\label{c+_def}
c_+(n,N)&=&\frac{N-n}{N}\frac{n\phi_A}{n\phi_A+(N-n-1)\phi_B}\ ,\\
\label{C0_def}
c_0(n,N)&=&\frac{n}{N}\frac{(n-1)\phi_A}{(n-1)\phi_A+(N-n)\phi_B}
+\frac{N-n}{N}\frac{(N-n-1)\phi_B}{n\phi_A+(N-n-1)\phi_B}\ ,\\
\label{c-_def}
c_-(n,N)&=&\frac{n}{N}\frac{(N-n)\phi_B}{(n-1)\phi_A+(N-n)\phi_B}\ .
\end{eqnarray}

From that, we may easily write an equation for the evolution of $P$:
\begin{eqnarray}
\nonumber
P(t+\Delta t,n,N)&=&c_+(n-1,N)P(t,n-1,N)+c_0(n,N)P(t,n,N)\\
\label{evol_disc}
&&+c_-(n+1,N)P(t,n+1,N)\ .
\end{eqnarray}
After imposing the boundary conditions $P(t,-1,N)=P(t,N+1,N)=0$, 
$\forall t\ge 0$, we conclude that the previous recursion is valid
for $t\ge 0$ and $n=0,1,\cdots, N$.

Originally, the Moran process was defined with 
a frequency-independent fitness, i.e., $\phi_{A,B}$ were
independent of the particular composition of the population.
We consider, however, the frequency-dependent case, and
we obtain the results for frequency-independent populations as a
special case.

Now, we obtain the fitness. For that, we first
consider a two players game, with
pay-off matrix given by: 
\begin{center}
\begin{tabular}{c|cc}
&I&II\\
\hline
I&$A$&$B$\\
II&$C$&$D$
\end{tabular}\ ,
\end{center}
where I and II are two pure strategies and $A,B,C,D>0$.
We call an $E_q$-strategist an
individual that plays I with probability $q$ and
II with probability $1-q$.

We assume that the two types play two (possibly) different
  strategies, $E_{q_1}$ and $E_{q_2}$.
The pay-off matrix is then given by
\begin{center}
\begin{tabular}{c|cc}
&$E_{q_1}$&$E_{q_2}$\\
\hline
$E_{q_1}$&$\widetilde A$&$\widetilde B$\\
$E_{q_2}$&$\widetilde C$&$\widetilde D$
\end{tabular}\ ,
\end{center}
where
\begin{eqnarray}
\label{tilA}
\widetilde A&\bydef&q_1^2A+q_1(1-q_1)(B+C)+(1-q_1)^2D\ ,\\
\label{tilB}
\widetilde B&\bydef&q_1q_2A+q_1(1-q_2)B+(1-q_1)q_2C+(1-q_1)(1-q_2)D\ ,\\
\label{tilC}
\widetilde C&\bydef&q_1q_2A+(1-q_1)q_2B+q_1(1-q_2)C+(1-q_1)(1-q_2)D\ ,\\
\label{tilD}
\widetilde D&\bydef&q_2^2A+q_2(1-q_2)(B+C)+(1-q_2)^2D\ .
\end{eqnarray}

For simplicity, we consider in this section only
pure strategists, i.e., $E_1$- and $E_0$-strategists
for type $\mathbb{A}$ and type $\mathbb{B}$ individuals respectively.
The general case follows easily from the results
in this section replacing $(A,B,C,D)$ by
$(\widetilde A,\widetilde B,\widetilde C,\widetilde D)$. 

We identify fitnesses and pay-offs, and then we have that
 the fitnesses for I- and II-strategists, for a population with
$n$ I-strategists, are given by 
\begin{eqnarray}
\label{phiA}
&&\phi_A=\frac{n-1}{N-1}A+\frac{N-n}{N-1}B\ ,\ \ n=1,\cdots,N\ ,\\
\label{phiB}
&&\phi_B=\frac{n}{N-1}C+\frac{N-n-1}{N-1}D\ ,\ \ n=0,\cdots,N-1\ .
\end{eqnarray}

Then, the evolution iteration is given by
Equation~(\ref{evol_disc}) with transition 
coefficients~(\ref{c+_def})--(\ref{c-_def}) and (\ref{phiA})--(\ref{phiB}).

\subsection{The discrete dynamics}

A natural question is what are the steady states of the iteration
defined by the Moran process. Here we show that the discrete model
cannot have a non-pure equilibrium.

Let us define the relative fitness as
\[
\rho_N(n)=\frac{\phi_A(n)}{\phi_{B}(n)}=\frac{(A-B)n+BN-A}{(C-D)n+(N-1)D}>0.
\]
Also, let
\[
f_N(n)=\frac{n}{N}\left(\frac{N-n}{N}\right)\quad\text{and}\quad
g_N(n,\rho)=\frac{N-1+(\rho-1)n}{N}.
\]
Then it is a straightforward computation to verify that 
\[
c_+(n,N)=\frac{f_N(n)\rho_N(n)}{g_N(n,\rho_N(n))},\quad
c_-(n,N)=\frac{f_N(n)}{g_N(n-1,\rho_N(n))}
\]
and
\[
c_0(n,N)=1-f_N(n)\left(\frac{\rho_N(n)}{g_N(n,\rho_N(n))}+\frac{1}{g_N(n-1,\rho_N(n))}\right). 
\]
Let $\mat$ 
be the iteration matrix of (\ref{evol_disc}). Then $\mat$ is a $N+1\times
N+1$, tridiagonal matrix, with entries given by
\[
\mat_{ii}=c_0(i,N),\quad i=0,\ldots, N,
\]
\[
\mat_{(i+1)i}=c_{+}(i,N)\quad\text{and}\quad
\mat_{i(i+1)}=c_{-}(i+1,N),\quad  i=0,\ldots, N-1.
\]
From this, and the fact that $\rho_N(n)>0$, it is easy to see that $\mat$ is
a nonnegative matrix. Since $c_{0}(n,N)+c_{+}(n,N)+c_{-}(n,N)=1$,
$\mat$ is column stochastic. 

The answer to question raised in the beginning of this section
is given by the following result:
\begin{prop}
\label{Akto0}
Let $\mat$ be as above and let
$\mathbf{P(t)}=(P(t,0),P(t,1),\ldots,P(t,N))^\dagger$.  Then
\begin{enumerate}
\item
\[
\lim_{k\to\infty}\mat^k=\begin{pmatrix}
1&1-F_1&\ldots&1-F_{N-1}&0\\
0&0&\ldots&\ldots&0\\
\vdots&\vdots&\vdots&\vdots&\vdots\\
0&0&\ldots&\ldots&0\\
0&F_1&\ldots&F_{N-1}&1
\end{pmatrix},
\label{our_prop}
\]
where the $F_n$ satisfy
\begin{align}
&F_{n}=c_{+}(n,N)F_{n+1}+c_{-}(n,N)F_{n-1}+c_0(n,N)F_n, \nonumber\\
&F_0=0\quad\text{and}\quad F_N=1.
\label{bs:fixeq}
\end{align}
\item If $\mathbf{1}$ denotes the vector $(1,1,\ldots,1)^\dagger$,
  $\mathbf{F}=(F_0,F_1,\ldots,F_N)^\dagger$ and  if
  $\langle\cdot,\cdot,\rangle$ denotes the usual inner  product, then 
we have that
\[
\langle \mathbf{P}(t),\mathbf{1} \rangle =\langle
\mathbf{P}(0),\mathbf{1} \rangle 
\quad\text{and}\quad
\langle \mathbf{P(t)}, \mathbf{F} \rangle = \langle \mathbf{P(0)},
\mathbf{F} \rangle.
\]
In particular, the $l^1$-norm of a nonnegative initial condition is preserved.
\end{enumerate}
\end{prop}
\begin{proof}
For part 1, see the proof at appendix~\ref{ProofProp1}.

As for part 2, we first observe that, if a vector $\mathbf{V}$ satisfies
$\mat^\dagger\mathbf{V}=\mathbf{V}$, then we have that
\[
\langle \mathbf{P(t+\Delta t)}, \mathbf{V} \rangle=
\langle \mathbf{M}\mathbf{P(t)}, \mathbf{V}
\rangle=
\langle \mathbf{P(t)}, \mathbf{M}^\dagger\mathbf{V}
\rangle=
\langle \mathbf{P(t)}, \mathbf{V}
\rangle.
\]
Hence 
\[
\langle \mathbf{P}(t),\mathbf{V} \rangle =\langle
\mathbf{P}(0),\mathbf{V}\rangle.
\]
From the fact that $\mat$ is column stochastic, we easily conclude
that
\[
\mat^\dagger\mathbf{1}=\mathbf{1},
\]
and the first invariant follows. For the second invariant, we observe that
Equation (\ref{bs:fixeq}) can be written in matrix notation as
\[
\mathbf{M}^{\dagger}\mathbf{F}=\mathbf{F},
\]
which concludes the proof.
\end{proof}

\begin{rmk}
The two invariants described in part 2 of the proposition \ref{Akto0}
are the only invariants of the Moran process and play an important
role in the determination of the correct continuous solution. 
\end{rmk}

Thus, the equilibrium states must have their mass concentrated in the
extremes. The $F_n$ turns out to be the fixation probability of
I-strategists, when the process start with $n$ I-strategists.

From the definitions of $c_{*}(n,N)$, we see that $F_n$ satisfies:
\begin{equation}
\left\{
\begin{array}{rcl}
&&\rho_N(n)F_{n+1}-\left(\rho_N(n)+\frac{g_N(n,\rho_N(n))}{g_N(n-1,\rho_N(n))}\right)F_n
+\frac{g_N(n,\rho_N(n))}{g_N(n-1,\rho_N(n))}F_{n-1}=0\ ,\\
&&F_0=0\quad \text{and}\quad F_N=1.
\end{array}
\right.
\label{gen:fixeq}
\end{equation}

Equation (\ref{gen:fixeq}) can be solved by writing
\[
H(n)=\frac{g_N(n,\rho_N(n))}{\rho_N(n)g_N(n-1,\rho_N(n))}
\quad\text{and}\quad G_n=F_n-F_{n-1}
\]
Then, ignoring the boundary conditions for the moment, we have that
\[
G_{n+1}=H(n)G_n,
\]
with solution given by
\[
G_n=G_1\prod_{i=1}^{n-1}H(i).
\]
Since
\[
F_{n}-F_{n-1}=G_1\prod_{i=1}^{n-1}H(i),
\]
we obtain, after applying $F_N=1$ and $F_0=0$, that:
\begin{align}
F_n&=G_1\sum_{k=1}^{n}\prod_{i=1}^{k-1}H(i), \nonumber\\
G_1&=\left(\sum_{k=1}^{N}\prod_{i=1}^{k-1}H(i)\right)^{-1}.
\label{gen:fixsoln}
\end{align}
The expression given by (\ref{gen:fixsoln}) does not  appear to yield
a simple formula in the general case. However, compare the
  formulas found in Section \ref{moran}, where
we study the case when the relative fitness is constant with respect
to $n$. 

\begin{rmk}
The  coefficients obtained in the above analysis are for a Death/Birth
process. For a Birth/Death process, they are simpler and are given by
\begin{eqnarray*}
c_+(n,N)&=&\frac{f_N(n)\rho_N(n)}{\tilde{g}_N(n,\rho_N(n))}\ ,\\
c_-(n,N)&=&\frac{f_N(n)}{\tilde{g}_N(n,\rho_N(n))}\ ,\\
c_0(n,N)&=&1-\frac{f_N(n)}{\tilde{g}_N(n,\rho_N(n))}\left(1+\rho_N(n)\right)\,
\end{eqnarray*}
where
\[
\tilde{g}_N(n,\rho)=\frac{N+(\rho-1)n}{N}.
\]
Also, in this case $H(n)$ simplifies to
\[
H(n)=\frac{1}{\rho_N(n)}.
\]
\end{rmk}

\subsection{Numerical Results}

We numerically computed the $\mat^{10000}$, for $N=20$ and various relative
fitnesses. The entries predicted to be zero by Proposition \ref{Akto0}
were found to have magnitude less than $10^{-50}$. 

\begin{figure}
\begin{center}
\epsfig{file=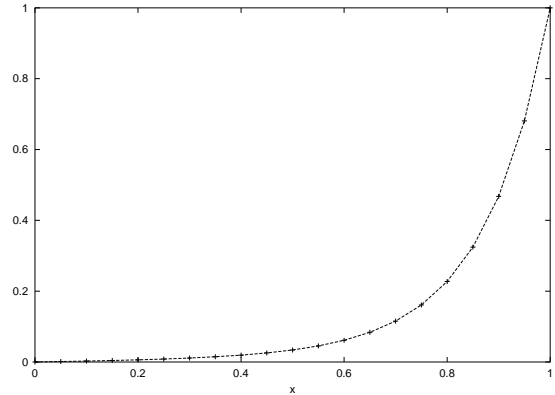,width=0.5\textwidth}
\end{center}
\caption{Fixation probabilities for $N=20$, $A=2$, $B=1$, $C=3$ and $D=1$.
The points are taken from  $\mat^{10000}$, while the lines are
obtained by numerically solving   (\ref{gen:fixsoln}).
\label{fix2fig}} 
\end{figure}

\begin{figure}
\begin{center}
\epsfig{file=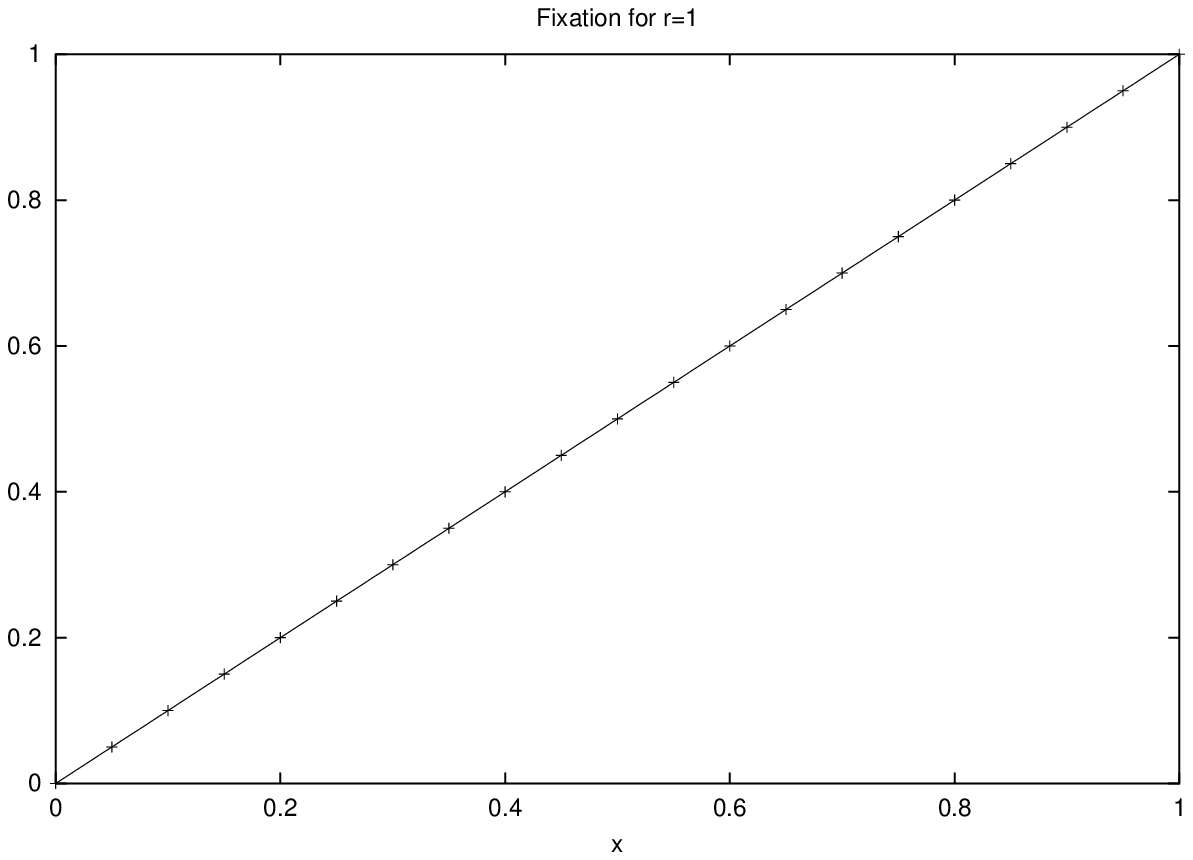,width=0.4\textwidth}
\hspace{0.1\textwidth}
\epsfig{file=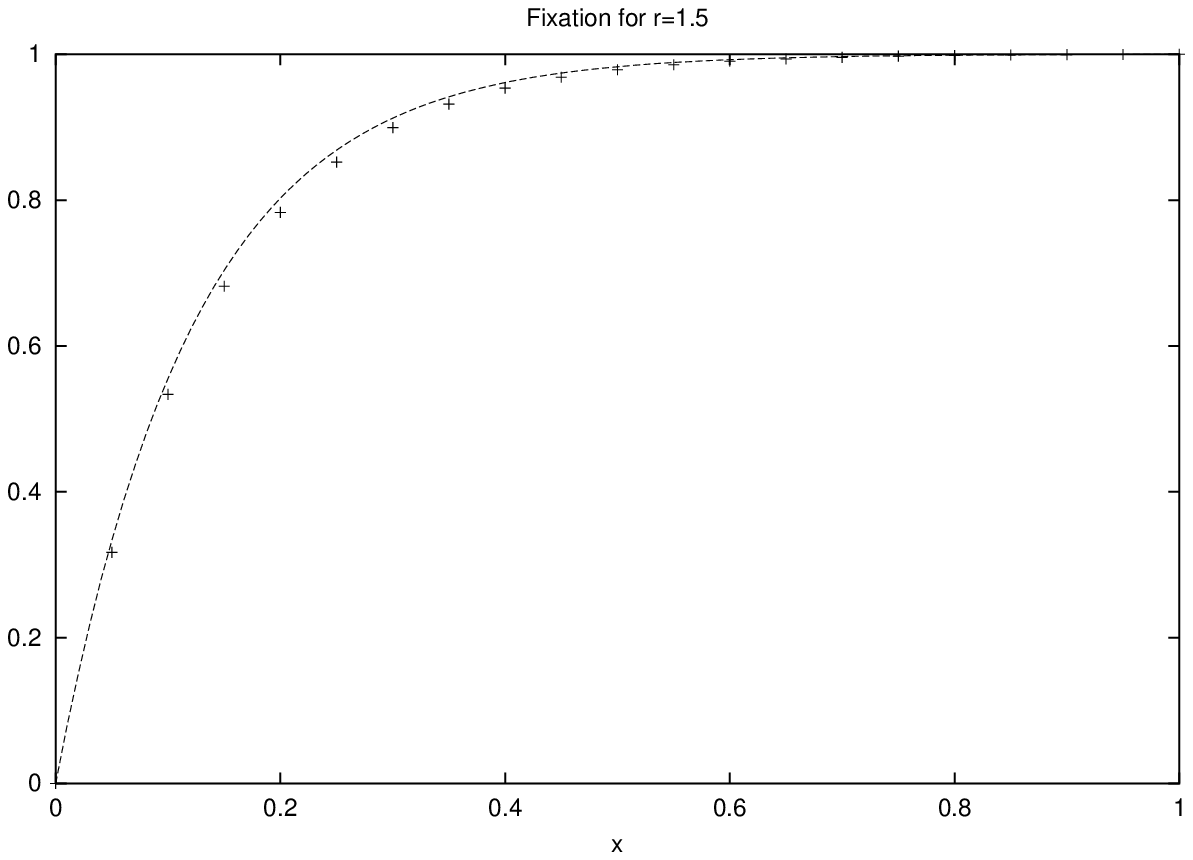,width=0.4\textwidth}
\vspace{.1cm}
\text{\hfil (a)\hspace{.45\textwidth} (b)\hfil}
\end{center}
\caption{Fixation probabilities for constant $r=C/A=D/B$, 
when $N=20$ computed
  from $\mat^{10000}$ together with the analytical fixation plotted as
  continuous functions  of $n/N$; (a) $r=1$ (b) $r=1.5$.
\label{fixfig}}
\end{figure}

Also, from these calculations, we extracted the fixation probabilities and
compared them with the ones obtained by evaluating (\ref{gen:fixsoln})
numerically. The result for a specific choice of fitness is displayed
in Figure \ref{fix2fig}.
For the case of frequency independent fitness, we can obtain explicit
formulas for the fixation probability--- see Section \ref{moran}
--- and we also 
compare with the fixation probabilities extracted from
$\mat^{10000}$ in  Figure~\ref{fixfig}.

\section{The thermodynamical limit}
\label{sec:thermo}

The aim of this section is to derive a continuous
approximation, i.e., a PDE model for the discrete process described in
the previous section.

We define the probability density that at time $t$ we
have a fraction $x\in[0,1]$ of type $\mathbb{A}$ individuals
\begin{equation*}
\mathcal{P}(t,x,N)\bydef \frac{P(t,xN,N)}{\frac{1}{N}}=NP(t,xN,N)\;,
\quad\text{with}\quad  x=\frac{n}{N},\quad n=0,1,2,\ldots,N.
\end{equation*}
Furthermore, 
we assume that in the limit $N\to\infty$, $\mathcal{P}(t,x,N)$ converges in
  some sense to a function $p(t,x)$ which is sufficiently smooth so that 
\begin{equation}
p\left(t,x\pm\frac{1}{N}\right)=p(t,x)\pm\frac{1}{N}\partial_xp(t,x)+
\frac{1}{2N^2}\partial_x^2p(t,x)+\O(N^{-3})
\end{equation}
and  re-write equation~(\ref{evol_disc}) to second order in
$N^{-1}$ as
\begin{eqnarray}
\label{EqForp}
&&p(t+\Delta t,x)-p(t,x)=\frac{1}{N}\left
[\left(c_+^{(1)}+c_0^{(1)}+c_-^{(1)}\right)p
-\left(c_+^{(0)}-c_-^{(0)}\right)\partial_xp\right]\\
\nonumber
&&\qquad
+\frac{1}{N^2}\left[\frac{1}{2}\left(c_+^{(2)}+c_0^{(2)}+c_-^{(2)}\right)p-
\left(c_+^{(1)}-c_-^{(1)}\right)\partial_xp
+\frac{1}{2}\left(c_+^{(0)}+c_-^{(0)}\right)\partial_x^2p\right]\\
\nonumber
&&\qquad
+\O\left(\frac{1}{N^3}\right)\ ,
\end{eqnarray}
where $c_*^{(i)}=c_*^{(i)}(x)$, $*=+,0,-$, $i=0,1,2$, are defined by
\begin{eqnarray}
\label{exp_c+}
c_+\left(\left(x-\frac{1}{N}\right)N,N\right)&=&c_+(n-1,N)
=c_+^{(0)}+\frac{1}{N}c_+^{(1)}+\frac{1}{2N^2}c_+^{(2)}\ ,\\
\label{exp_c0}
c_0(xN,N)&=&c_0(n,N)=c_0^{(0)}+\frac{1}{N}c_0^{(1)}+
\frac{1}{2N^2}c_0^{(2)}\ ,\\
\label{exp_c-}
c_-\left(\left(x+\frac{1}{N}\right)N,N\right)&=&c_-(n-1,N)=
c_-^{(0)}+\frac{1}{N}c_-^{(1)}+\frac{1}{2N^2}c_-^{(2)}\ .
\end{eqnarray}
Then
\begin{eqnarray}
\label{c1soma}
c_+^{(1)}+c_0^{(1)}+c_-^{(1)}&=&
\left(Ax^2+D(1-x)^2+(B+C)x(1-x)\right)^{-2}\cdot
\bigl[A(A-C)x^4\\
\nonumber
&&+\left(B(B-D)+C(C-A)+2C(B-D)\right)x^2(1-x)^2\\
\nonumber
&&+D(D-B)(1-x)^4+2x(1-x)\left(A(B-D)x^2-(A-C)D(1-x)\right)\bigr]\ ,\\
\label{c0menos}
c_+^{(0)}-c_-^{(0)}&=&\frac{x(1-x)\left(x(A-C)+(1-x)(B-D)\right)}
{Ax^2+D(1-x)^2+(B+C)x(1-x)}\ ,
\end{eqnarray}
If we impose that
\begin{eqnarray}
\label{Ntoinfty1}
&&\lim_{N\to\infty}(A,B,C,D)=(1,1,1,1)\ ,\\
\label{Ntoinfty2}
&&\lim_{N\to\infty}N(A-1,B-1,C-1,D-1)=(a,b,c,d)\ ,
\end{eqnarray}
we find
\begin{eqnarray*}
\lim_{N\to\infty}\left(c_+^{(2)}+c_0^{(2)}+c_-^{(2)}\right)&=&-4\ ,\\
\lim_{N\to\infty}\left(c_+^{(1)}-c_-^{(1)}\right)&=&-2+4x\ ,\\
\lim_{N\to\infty}\left(c_+^{(0)}+c_-^{(0)}\right)&=&2x(1-x)\ ,
\end{eqnarray*}
and, from~(\ref{c1soma}--\ref{c0menos}), we have
\begin{eqnarray*}
\lim_{N\to\infty}N\left(c_+^{(1)}+c_0^{(1)}+c_-^{(1)}\right)&=&
-3x^2(a-b-c+d)-2x(a-c-2(b-d))+(d-b)\ ,\\
\lim_{N\to\infty}N\left(c_+^{(0)}-c_-^{(0)}\right)&=&
x(1-x)(x(a-c)+(1-x)(b-d))\ .
\end{eqnarray*}
Finally, we divide Equation~(\ref{EqForp}) by $\Delta t=N^{-2}$ 
(diffusive scaling),  
and take the limit $N\to\infty$, to obtain 
\begin{eqnarray*}
\partial_tp&=&\left[3x^2(a-b-c+d)-2x(a-c-2(b-d))-(b-d)\right]p\\
&&-x(1-x)(x(a-c)+(1-x)(b-d))\partial_xp\\
&&+(-2)p+2(1-2x)\partial_xp+x(1-x)\partial_x^2p
\end{eqnarray*}
i.e.,
\begin{equation}
\label{replicator_pde_pure}
\partial_t p=\partial_x^2\left[x(1-x)p\right]-
\partial_x\left[x(1-x)(x\alpha+(1-x)\beta)p\right]\ .
\end{equation}
where $\alpha=a-c$ and $\beta=b-d$.
We also define $\eta=\alpha-\beta$.

Supplementing Equation (\ref{replicator_pde_pure}) we have the
following conservation laws:
\[
\frac{\d}{\d t}\int_0^1p(t,x)\d x = 0\quad\text{and}\quad
\frac{\d}{\d t}\int_0^1\psi(x)p(t,x)\d x = 0,
\]
where $\psi(x)$ is given in Theorem~\ref{conv_fi}.

\begin{rmk}\label{scalling_rmk}
It is important to stress that if we do not impose conditions
(\ref{Ntoinfty1})--(\ref{Ntoinfty2}), there are another possible
scalings. More precisely, if  (\ref{Ntoinfty1}) still holds but
(\ref{Ntoinfty2}) is replaced by 
\[
\lim_{N\to\infty}N^\nu(A-1,B-1,C-1,D-1)=(a,b,c,d),\qquad  0<\nu<1,
\]
then another possible scaling is given by taking
$\Delta t=(1/N)^{1+\nu}$ and, in this case, we obtain
(\ref{replicator_pde_pure}) without the diffusion term. This equation
is discussed in Section~\ref{diffusionless}. Moreover, if we drop
(\ref{Ntoinfty1})--(\ref{Ntoinfty2}), and only require that the
payoffs have a finite limit when $N$ goes to infinity, then yet
another scaling is given by $\Delta t=1/N$ and, in this case, the
equation for the probability density is given by
\begin{equation}\label{drift-independent}
\partial_t\bar p=\partial_x\left[
\frac{x(1-x)\left(x(A-C)+(1-x)(B-D)\right)}{x^2(A-B-C+D)+x(B+C-2D)+D}
\bar p\right]\ .
\end{equation}
We analyze this equation in Section~\ref{drift}.
\end{rmk}

Equation (\ref{replicator_pde_pure}) is not readily covered by the usual
theory of parabolic PDEs. However, the analysis can be extended to
obtain the following result:

\begin{thm}
\begin{enumerate}
\item For a given $p^0\in L^1([0,1])$, there exists a unique solution
  $p=p(t,x)$  to Equation~(\ref{our_equation}) of class 
$C^\infty\left(\mathbb{R}^+\times(0,1)\right)$ that satisfies $p(0,x)=p^0(x)$.
\item The solution can be written as
\[
p(t,x)=q(t,x)+a(t)\delta_0+b(t)\delta_1,
\]
where $q\in C^{\infty}(\mathbb{R}^{+}\times [0,1])$ satisfies
(\ref{replicator_pde_pure}) without boundary conditions, and we also have
\[
a(t)=\int_0^tq(s,0)\d s\quad\text{and}\quad
b(t)=\int_0^tq(s,1)\d s.
\]
In particular, we have that $p\in C^{\infty}(\mathbb{R}^{+}\times (0,1))$.
\item We also have that
\[
\lim_{t\to\infty}q(t,x)=0\text{ (uniformly)},\quad \lim_{t\to\infty}
a(t)=\pi_0[p^0]\quad \text{and}\quad
 \lim_{t\to\infty} b(t)=\pi_1[p^0],
\]
where  $\pi_0$ and $\pi_1$ are computed in Theorem~\ref{conv_fi}.
Note that this means that the solution solution will 'die out' in
the interior and only the Dirac masses in the extremities will survive.

\item Assume $p^0\in L^2([0,1])$ and let
  $J(t)=\int_0^1x(1-x)q^2(t,x)\d x$. Then, 
  we have that  
\[
J(t)\leq J(0)e^{-2\lambda_0t}, \lambda_0>0.
\]
\end{enumerate}
\label{pde_prop}
\end{thm}

See the proof at Appendix~\ref{app_b}.

\begin{figure}
\begin{center}
\epsfig{file=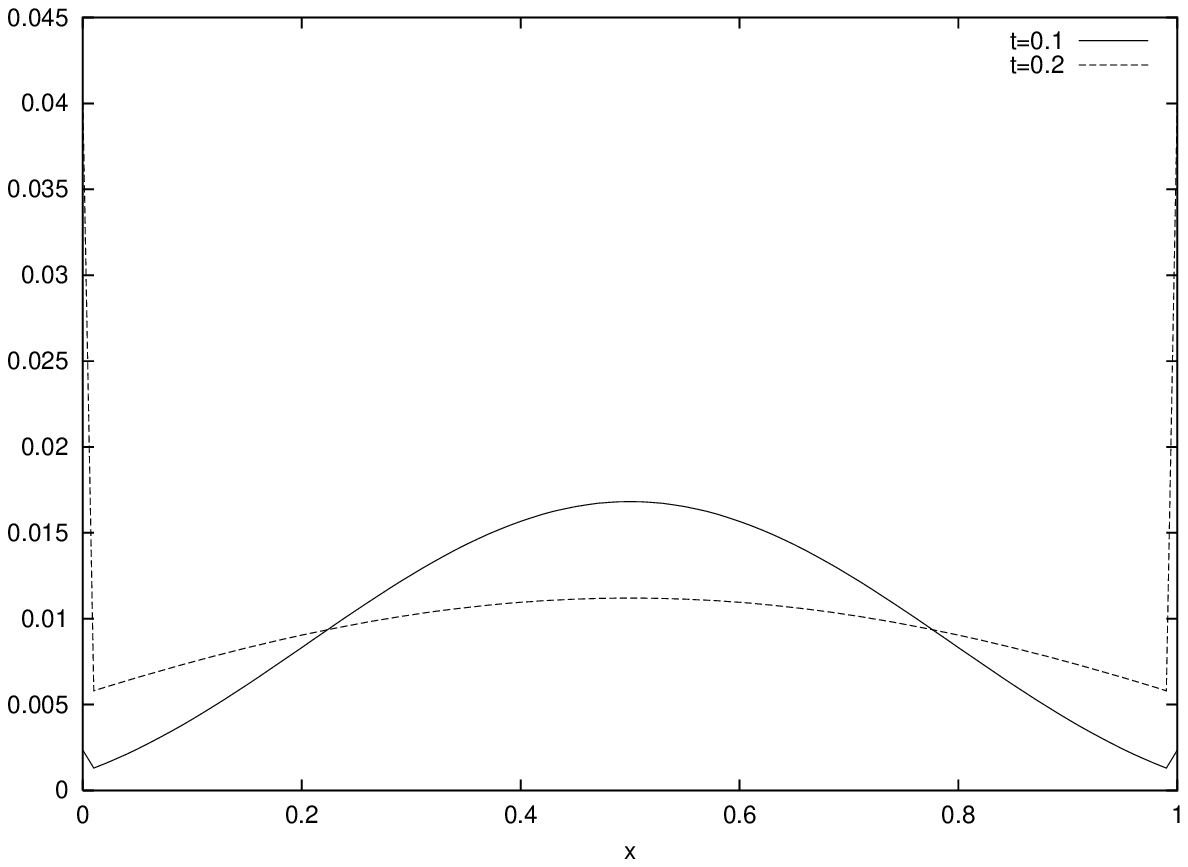,width=0.4\textwidth}
\hspace{0.2cm}
\epsfig{file=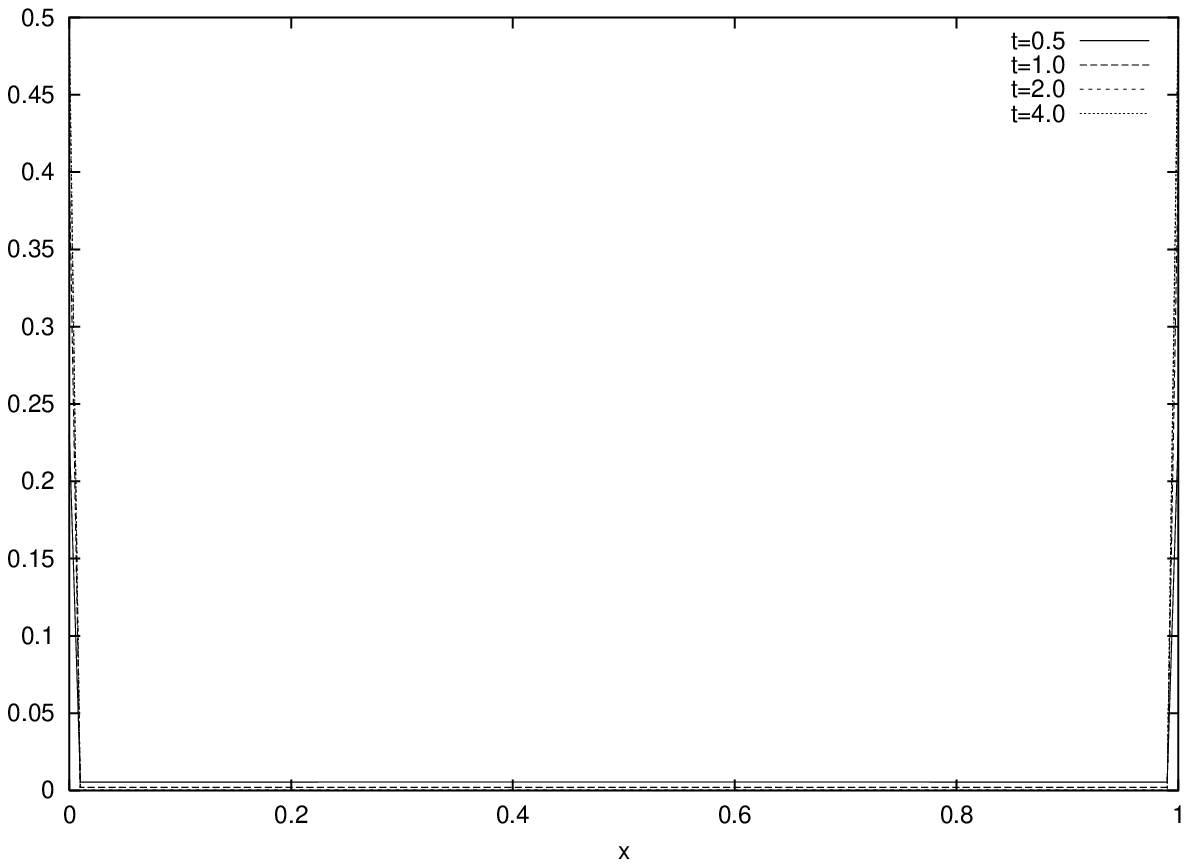,width=0.4\textwidth}
\end{center}
\caption{Solutions for $P(t,x)$ for various times, when
  $\beta=\alpha=0$. This is the pure diffusive constant fitness case.
Note the diffusion to the boundaries. The initial condition is given by
$p^0(x)=\delta_{1/2}(x)$.
\label{nfig}}
\end{figure}

\begin{figure}
\begin{center}
\epsfig{file=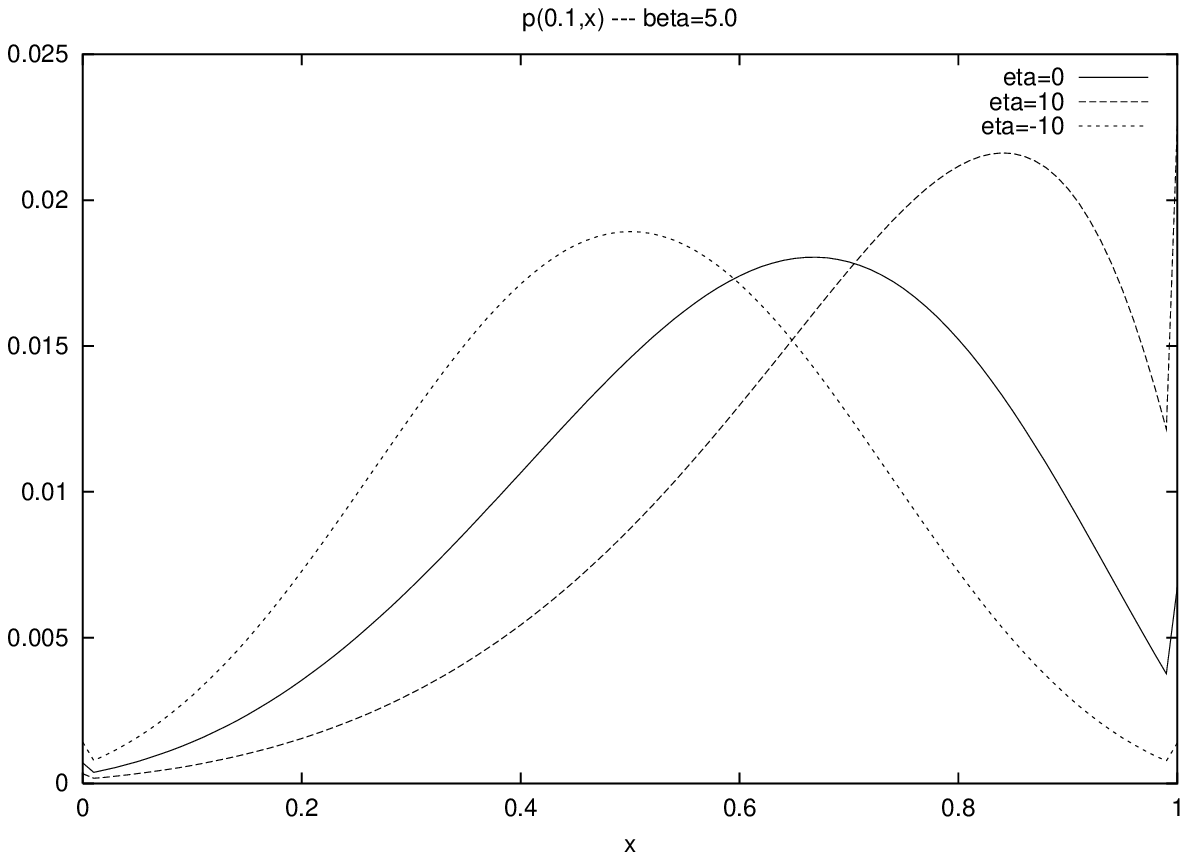,width=0.4\textwidth}
\hspace{0.2cm}
\epsfig{file=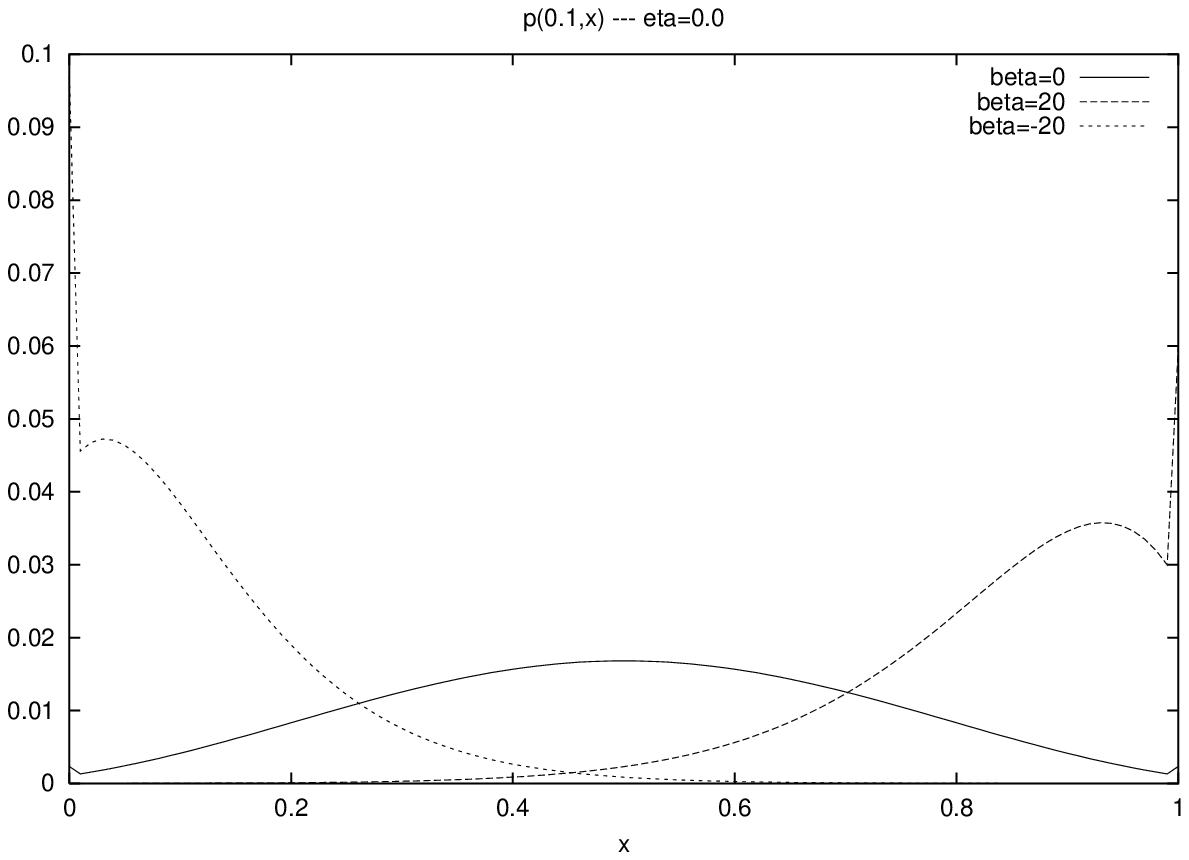,width=0.4\textwidth}
\end{center}
\caption{Solutions for $P(0.1,x)$ for various values $\beta$ and
  $\eta\bydef\alpha-\beta$. Here, the initial condition is 
the same as in Figure~\ref{nfig}. \label{nnfig}}
\end{figure}

\begin{figure}
\begin{center}
\epsfig{file=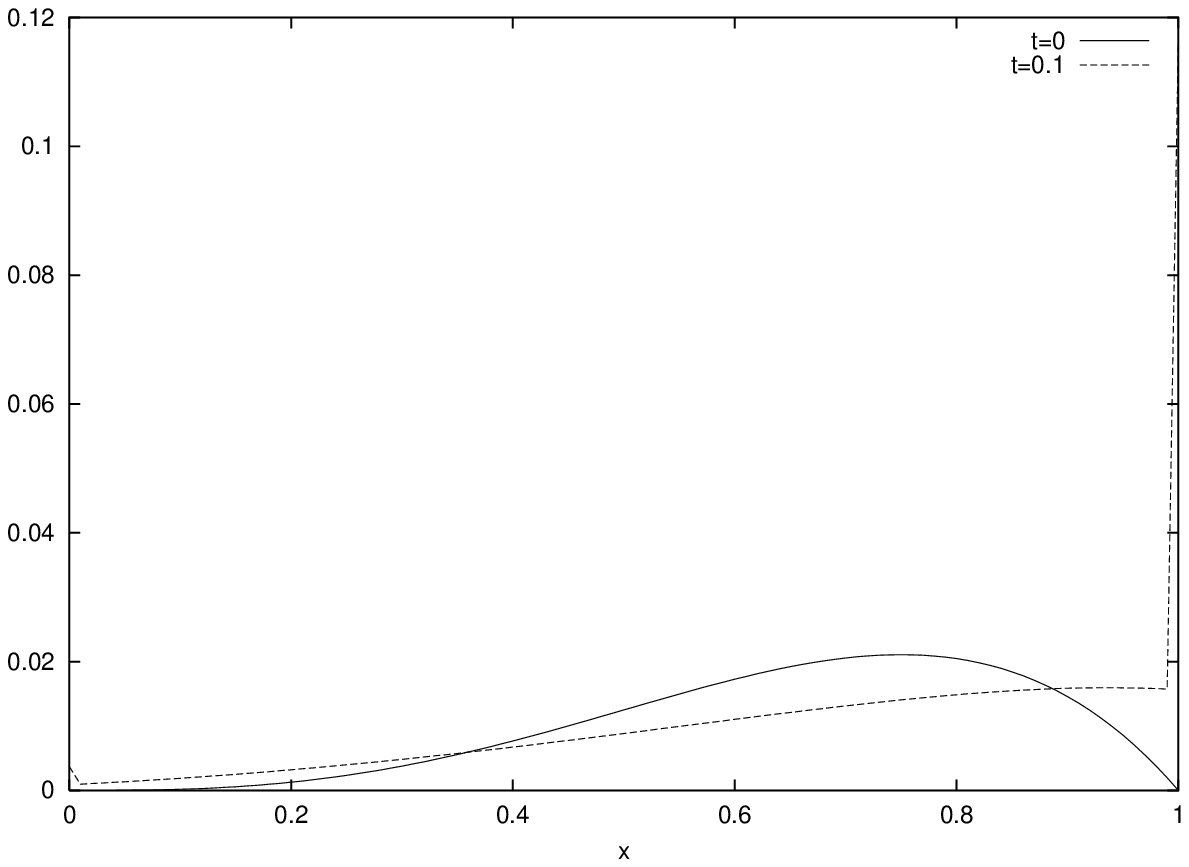,width=0.4\textwidth}
\hspace{0.2cm}
\epsfig{file=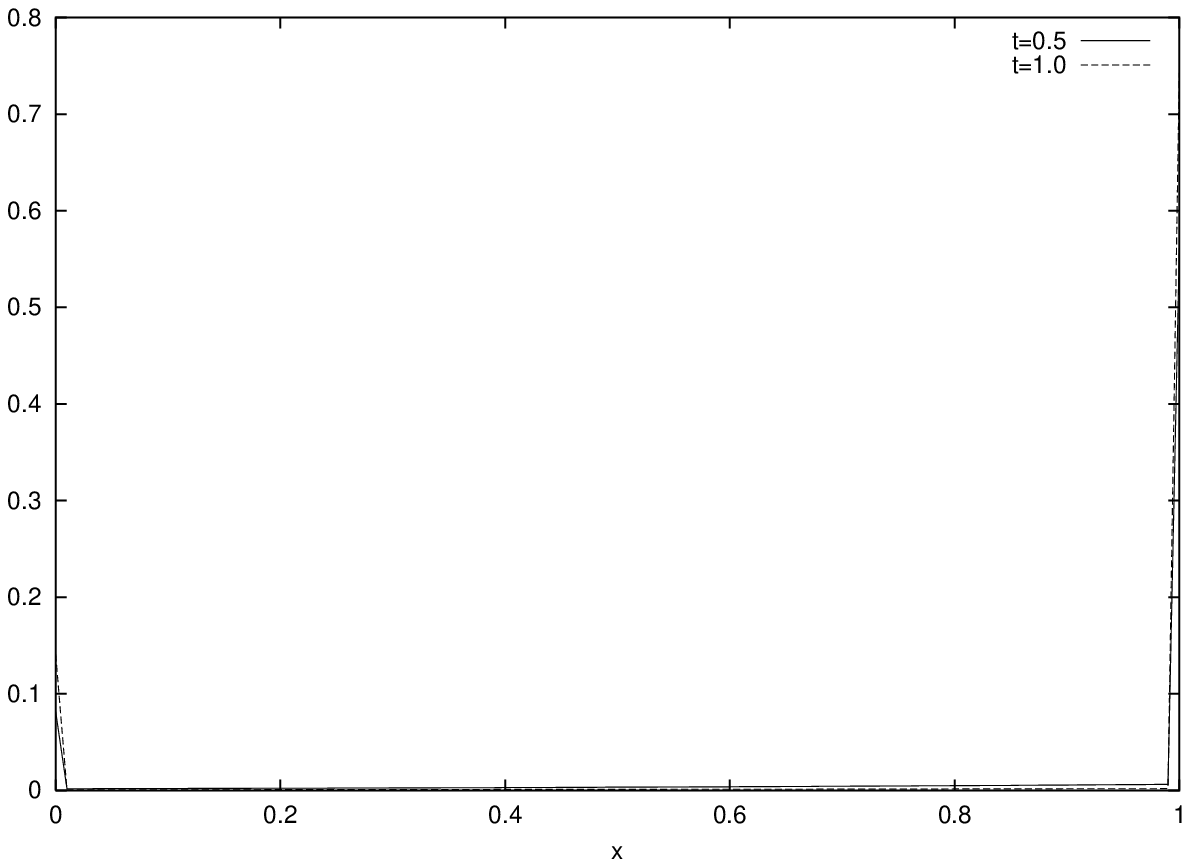,width=0.4\textwidth}
\end{center}
\caption{Solutions for $P(t,x)$ when $\beta=2$ and $\eta=0$ 
for various times. This is the case of some drift with constant
fitness. The initial condition is $p^0(x)=20x^3(1-x)$, which is
asymmetric with a peak at $x=3/4$. Notice that the form of the
initial condition together with the drift sign leads to
  a very rapid convergence to the equilibrium state.\label{mfig}}
\end{figure}

\begin{figure}
\begin{center}
\epsfig{file=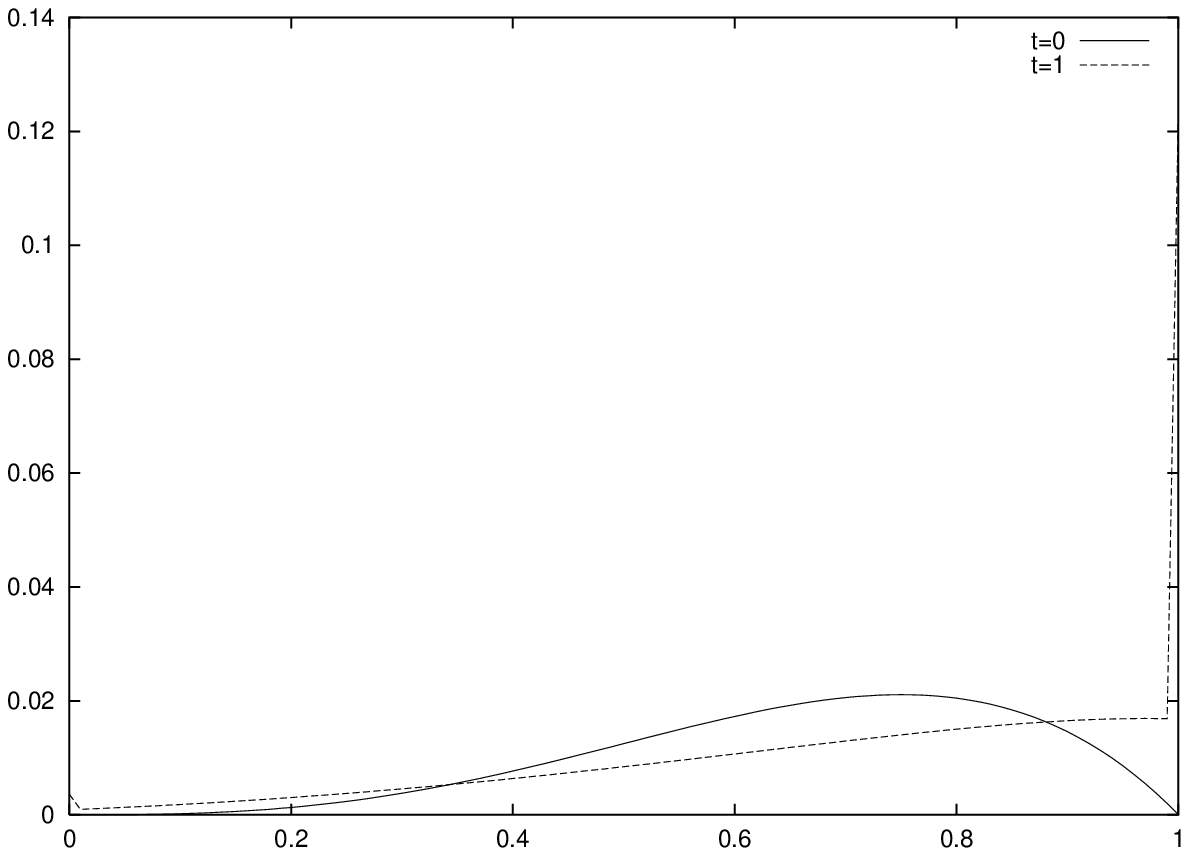,width=0.4\textwidth}
\hspace{0.2cm}
\epsfig{file=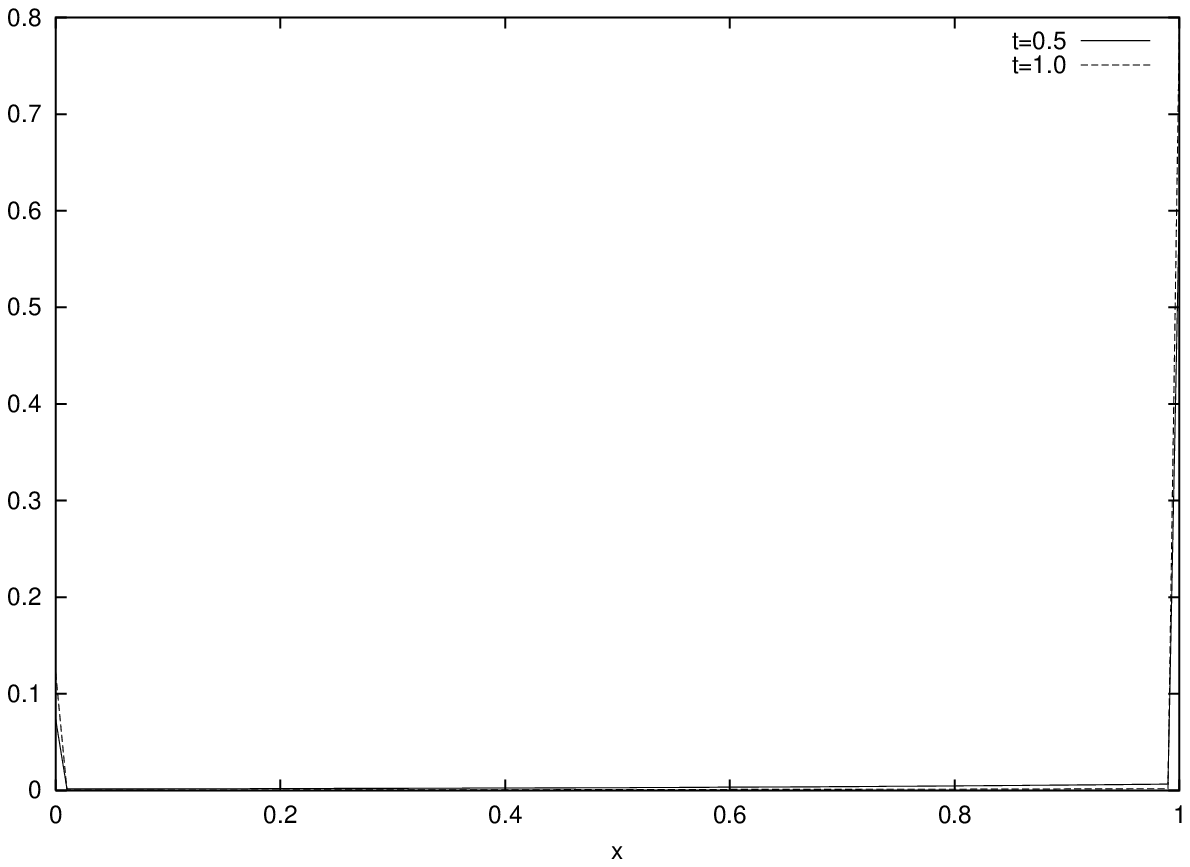,width=0.4\textwidth}
\end{center}
\caption{Solutions for $P(t,x)$ for various times, when $\beta=2$ and
  $\eta=1$. The initial condition is the same as in in figure~\ref{mfig}.
Notice that there is little difference from the  computation with
  $\eta=0$ thanks to the form of the initial   condition and
  to the  order one size of the parameters. \label{mmfig}} 
\end{figure}

\begin{figure}
\begin{center}
\epsfig{file=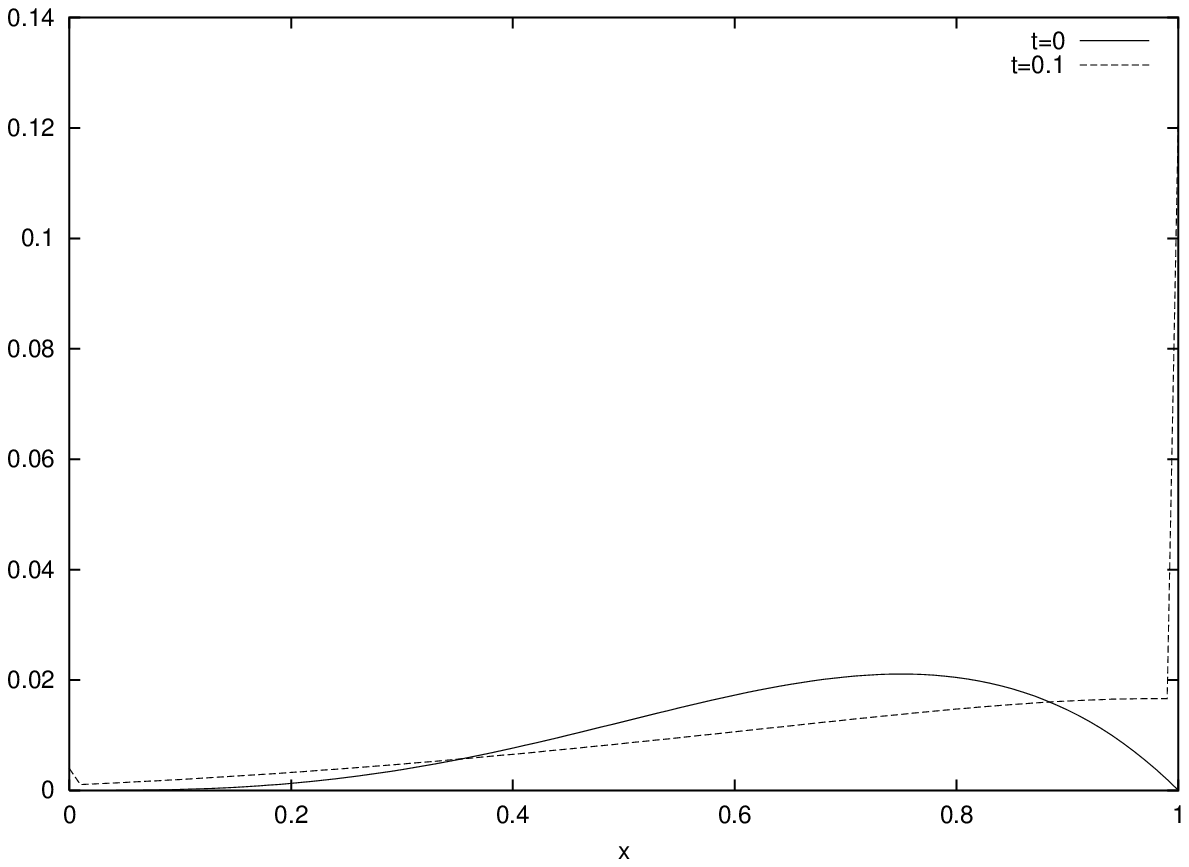,width=0.4\textwidth}
\hspace{0.2cm}
\epsfig{file=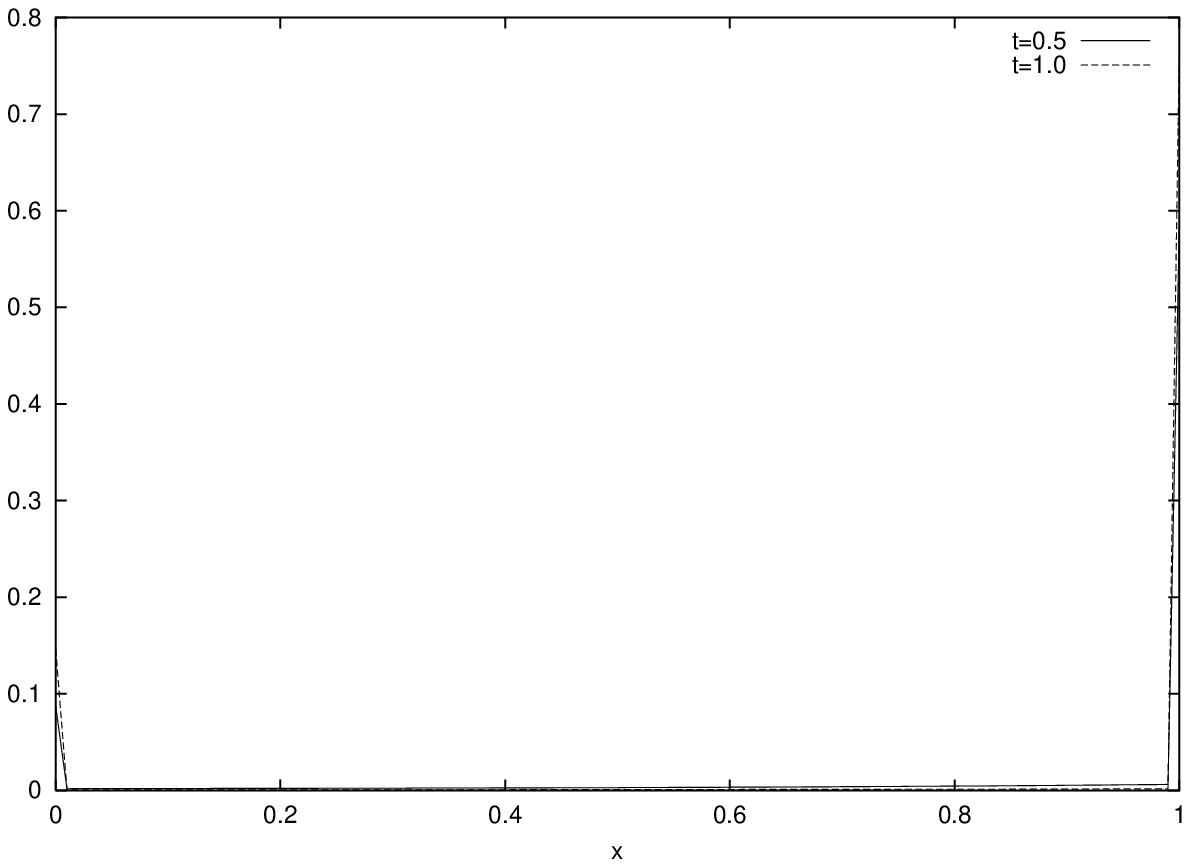,width=0.4\textwidth}
\end{center}
\caption{Same as Figure \ref{mmfig}, but with $\beta=1$ and
  $\eta=2$. Same remarks apply in this case. \label{mmmfig}}
\end{figure}

\begin{figure}
\begin{center}
\epsfig{file=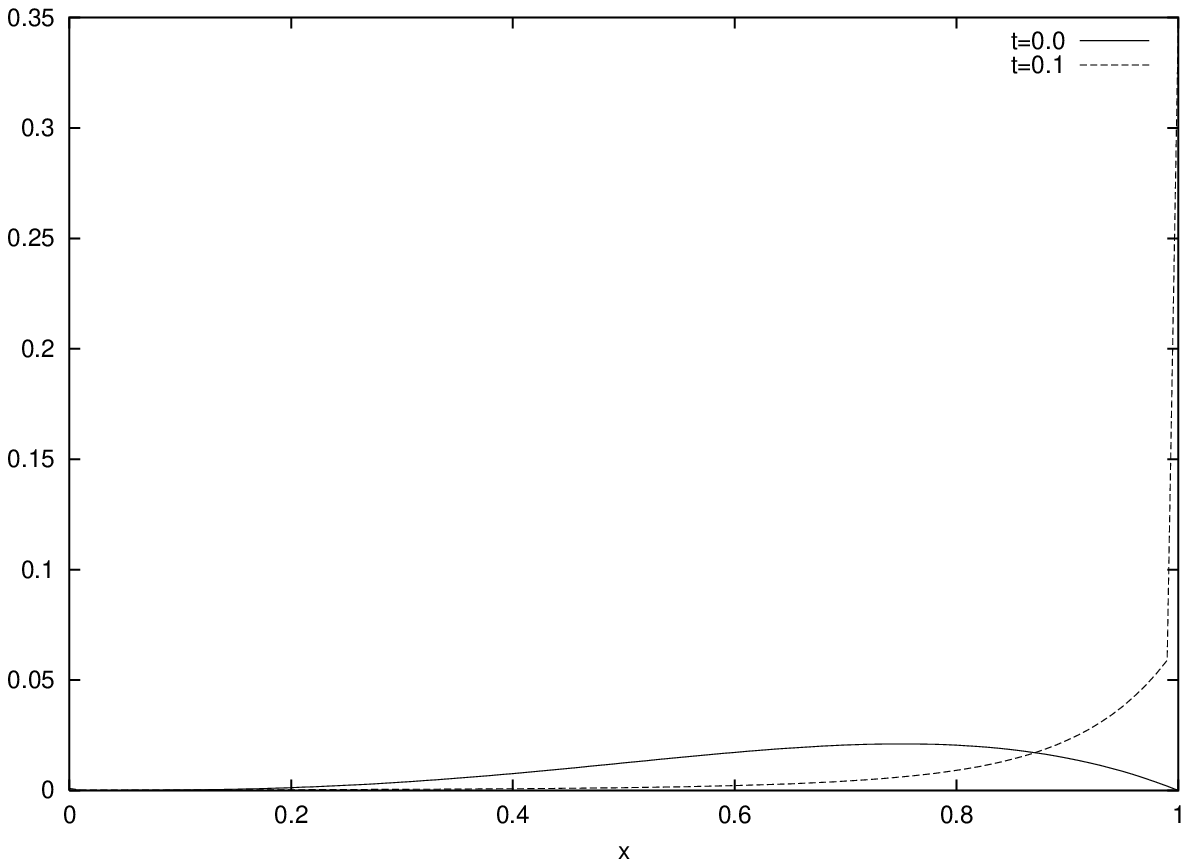,width=0.4\textwidth}
\hspace{0.2cm}
\epsfig{file=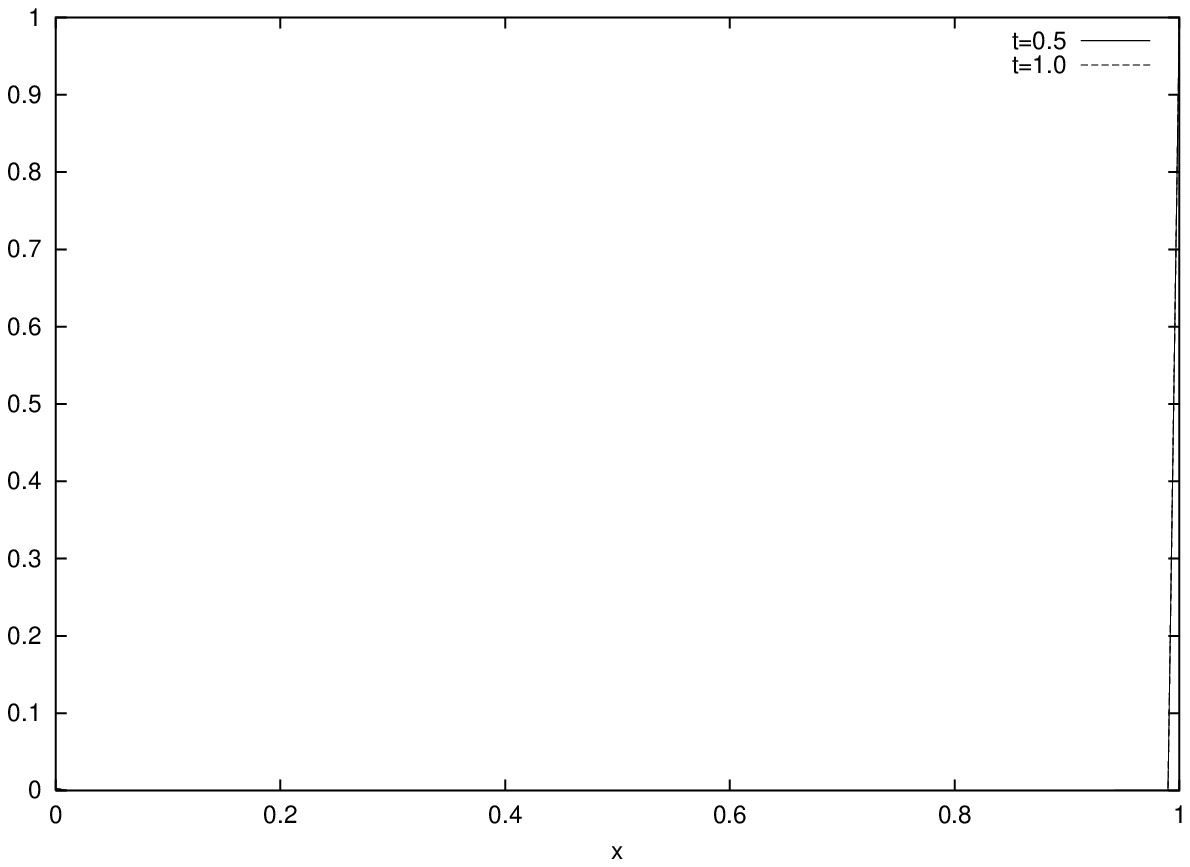,width=0.4\textwidth}
\end{center}
\caption{Solutions for $P(t,x)$ for various times when $\beta=10$ and
  $\eta=20$ with the same initial condition as in Figure~\ref{mfig}. 
The   convergence for the equilibrium state is very fast also in this
  case. \label{mmmmfig}}
\end{figure}

\begin{figure}
\begin{center}
\epsfig{file=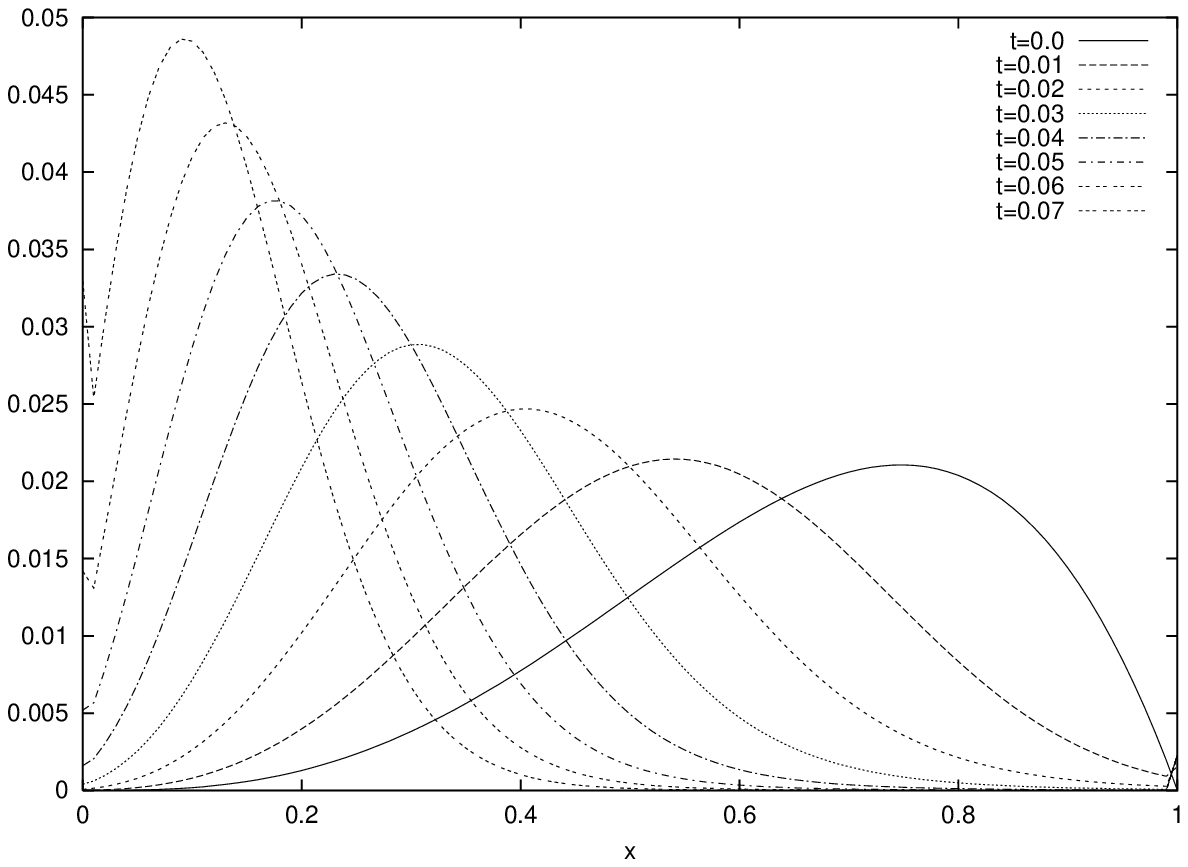,width=0.4\textwidth}
\hspace{0.2cm}
\epsfig{file=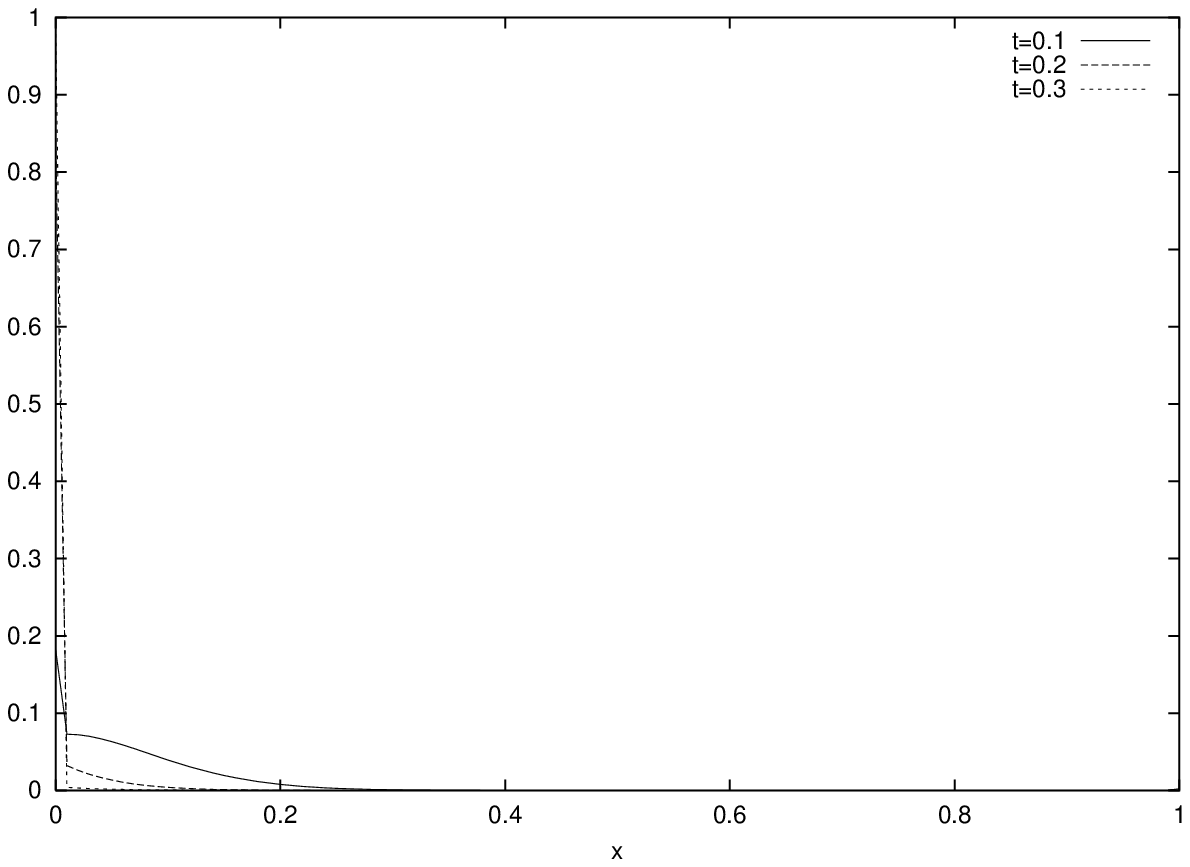,width=0.4\textwidth}
\end{center}
\caption{Solutions for $P(t,x)$ for various times when $\beta=-20$ and
  $\eta=-40$. In this case the drift forces the solution to
  accumulate in the opposite direction of the   initial large
  concentration. \label{mmmmmfig}} 
\end{figure}

\begin{figure}
\begin{center}
\epsfig{file=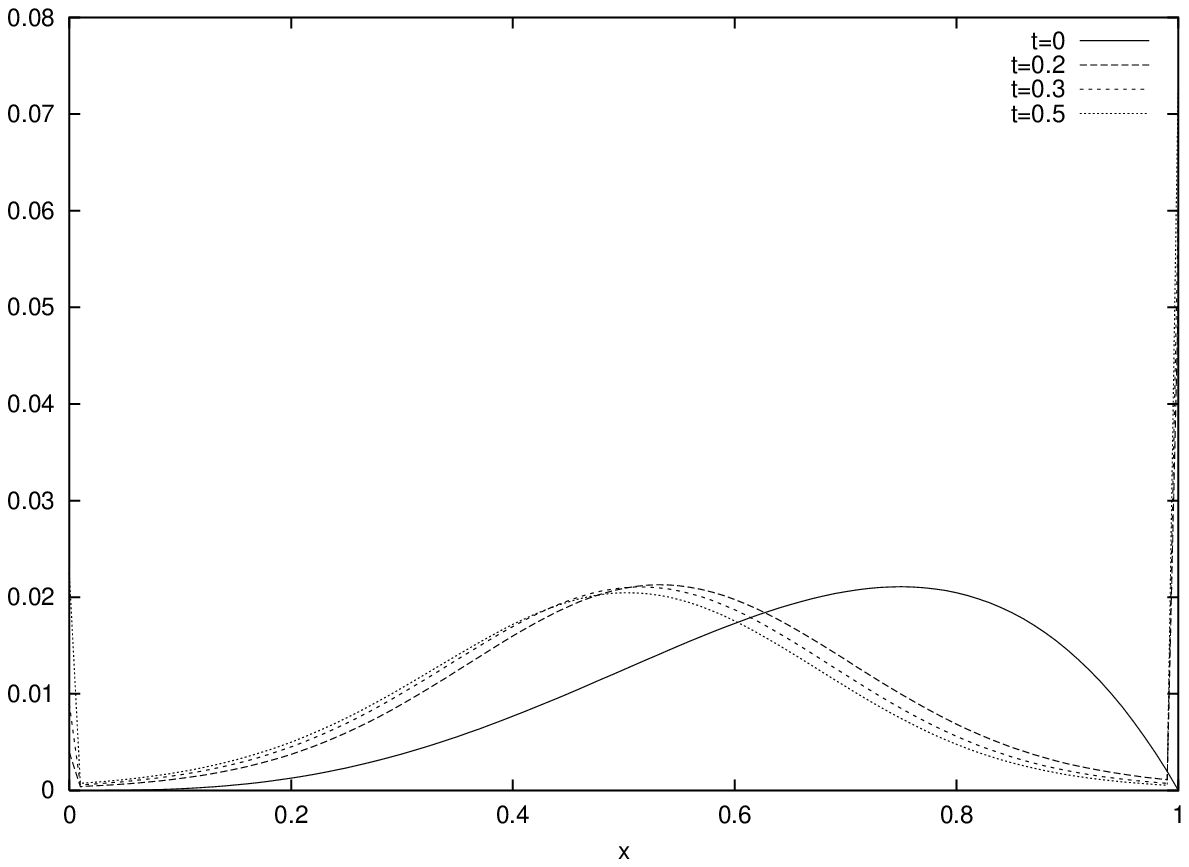,width=0.4\textwidth}
\hspace{0.2cm}
\epsfig{file=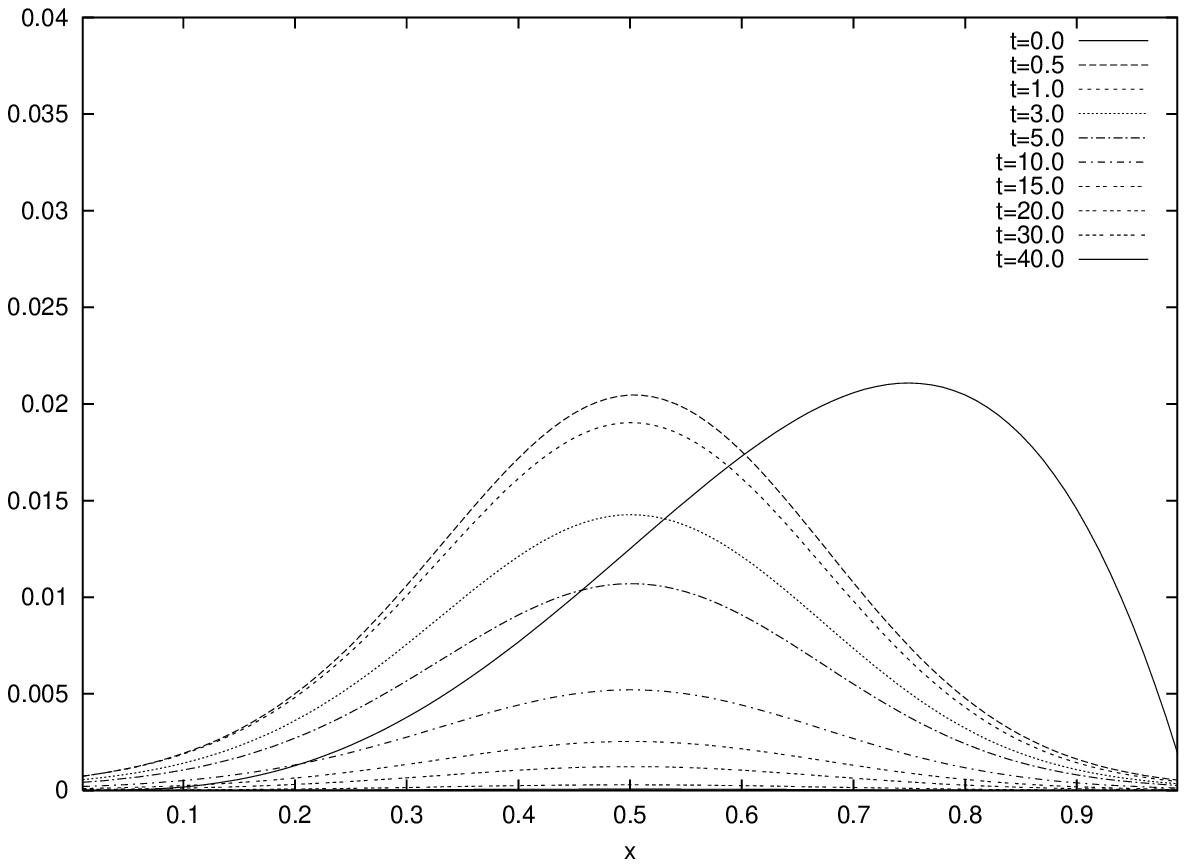,width=0.4\textwidth}
\end{center}
\caption{Solutions for $P(t,x)$ for various times when $\beta=20$ and
  $\eta=-40$. Here, the convective term vanishes at
  $x=1/2$. The effect is that, at first, the solution convected until
  the its peak reaches $x=1/2$. Then it essentially stays there, while
  diffusion enforces the transport to the boundaries. In the second
  figure, the very ends of the interval are omitted for better view of
  the behavior in interior. \label{mmmmmmfig}} 
\end{figure}

\begin{figure}
\begin{center}
\epsfig{file=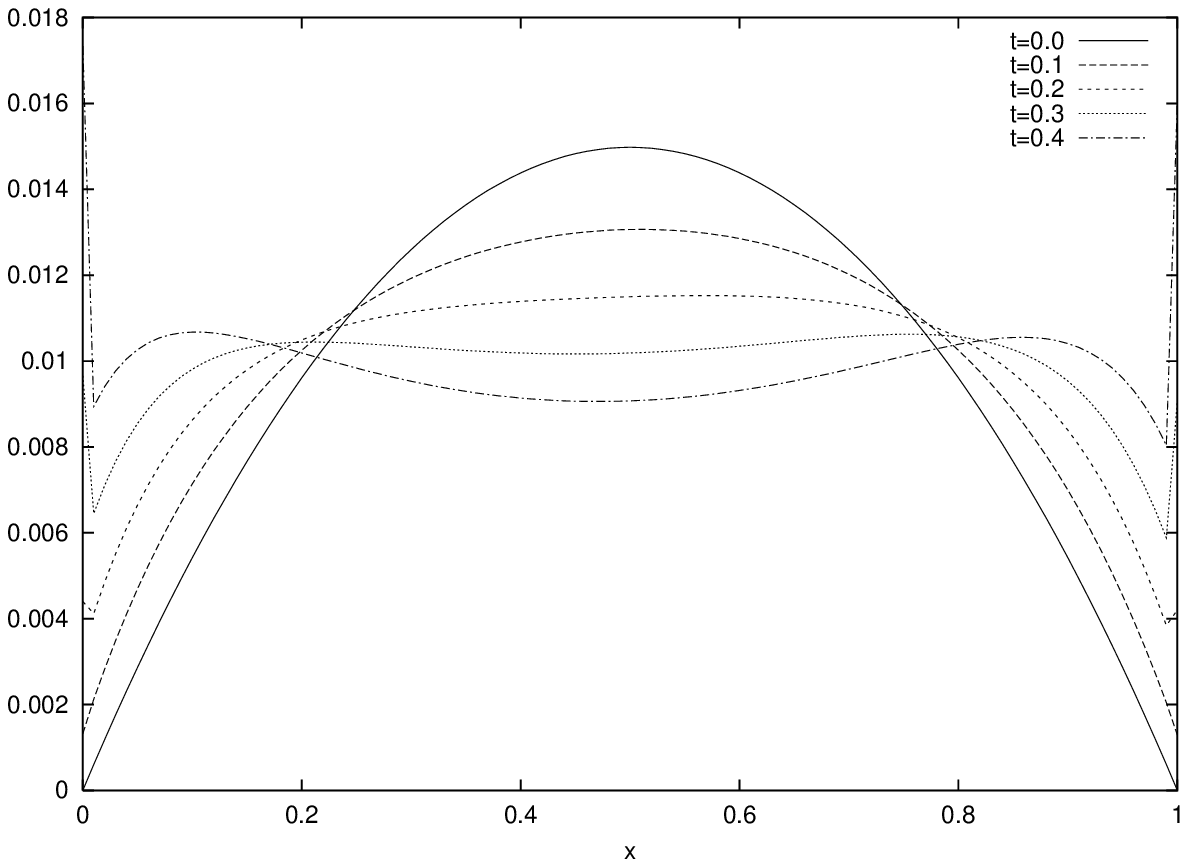,width=0.4\textwidth}
\hspace{0.2cm}
\epsfig{file=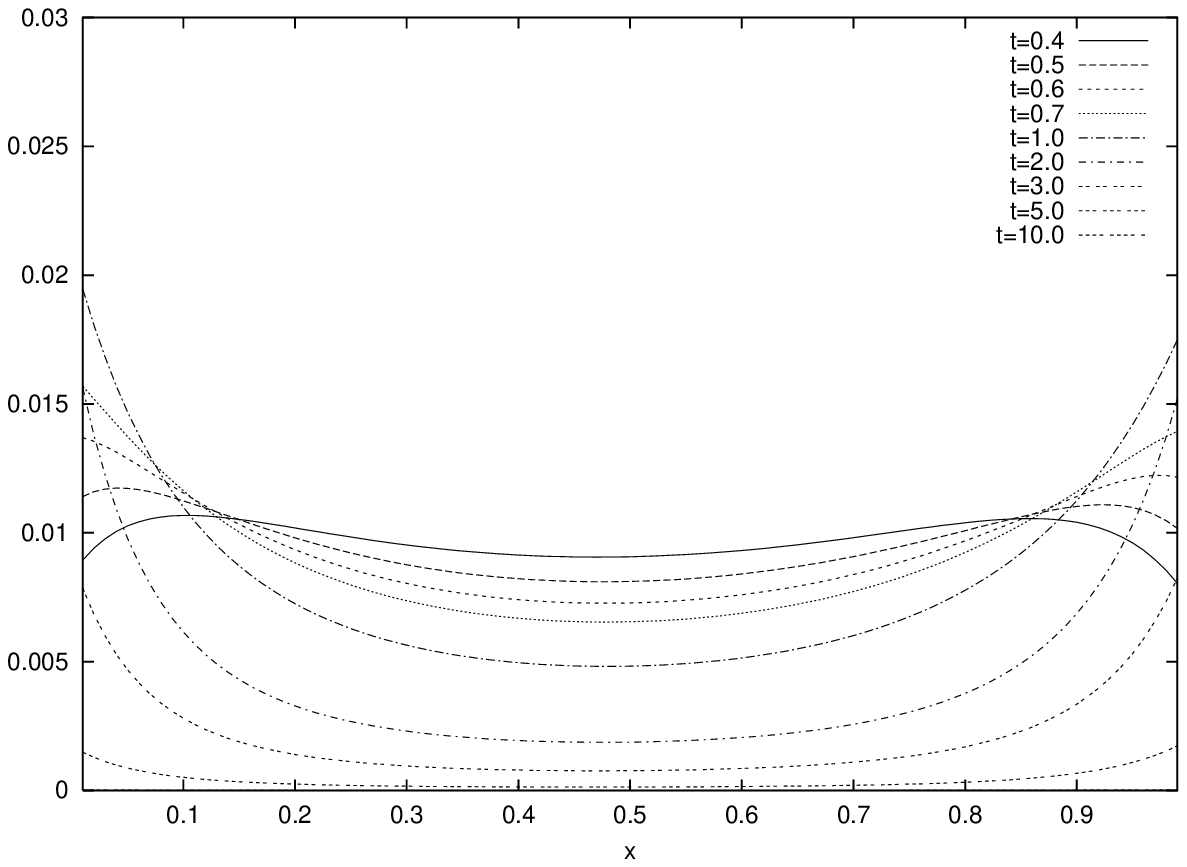,width=0.4\textwidth}
\end{center}
\caption{Solutions for $P(t,x)$ for various times, when $\beta=-20$ and
  $\eta=40$, with initial condition $p^0(x)=6x(1-x)$. In this case,
  the sign of $\eta$ drives the solution out of $x=1/2$ to the
  extremes.
The initial condition was chosen to be symmetric in
  this case to highlight this behavior. Also, as in the previous
  example, the second figure have the very ends of the interval
  omitted for a better view of the inner behavior.\label{mmmmmmmfig}}
\end{figure}

For completeness we show various  numerical simulations for computing
$p(t,x)$. Due to display convenience we plot $P(t,x)=(\Delta x)
p(t,x)$, instead of $p(t,x)$. See Figures~\ref{nfig}--\ref{mmmmmmmfig}.

We observe that $p'(t,x)=p(t,1-x)$ also satisfies
(\ref{replicator_pde_pure})  changing the
parameters $(\alpha,\beta)\to(-\beta,-\alpha)$.
Hence each computation actually yields solution for two set
of parameters, just by reflecting the solution around the axis $x=1/2$.


Notice that, for $\alpha-\beta,\beta\gg 1$,
we expect a behavior drift-dominated  for intermediate times. This
means that, if  $x^*=\beta/(\beta-\alpha)\not\in(0,1)$, then
the solution will be convected until it reaches one of the boundaries,
and then will diffuse to the steady state. Otherwise, 
depending on the sign of $\eta$
the solution will 
either first concentrate near $x^*$, and then diffuses to the
boundary, or depart from $x^*$ in both directions towards the ends.
Notice also, that the solutions are 
never smooth at the ends. Computations with different values of
$\alpha$ and $\beta$ produces qualitatively similar graphics.

Now, let us go back to the general case, i.e., for
$E_{q_1}$- and $E_{q_2}$-strategists, instead of only pure strategists. 
Then, in a straightforward way 
(see also Equations~(\ref{tilA})--(\ref{tilD})), 
we define
\begin{eqnarray*}
\tilde a&\bydef&q_1^2a+q_1(1-q_1)(b+c)+(1-q_1)^2d\ ,\\
\tilde b&\bydef&q_1q_2 a +q_1(1-q_2)b+(1-q_1)q_2c+(1-q_1)(1-q_2)d\ ,\\
\tilde c&\bydef&q_1q_2a +(1-q_1)q_2b+q_1(1-q_2)c+(1-q_1)(1-q_2)d\ ,\\
\tilde d&\bydef&q_2^2a+q_2(1-q_2)(b+c)+(1-q_2)^2d\ .
\end{eqnarray*}

Then, the equation for $p$, the fraction of $E_{q_1}$-strategists
in the population is given by
\begin{eqnarray}
\nonumber
\partial_t p&=&\partial_x^2\left(x(1-x)p\right)
-\partial_x\left(x(1-x)(x(\tilde a-\tilde c)+
(1-x)(\tilde b-\tilde d))p\right)\\
\label{replicator_pde}
&=&\partial_x^2\left(x(1-x)p\right)
-\partial_x\left(x(1-x)(x(\tilde\alpha+(1-x)\tilde\beta))p\right)\ ,
\end{eqnarray}
where $\tilde\alpha\bydef \tilde a-\tilde c
=(q_1-q_2)(q_1\alpha+(1-q_1)\beta)$ and 
$\tilde\beta\bydef\tilde b-\tilde d=(q_1-q_2)(q_2\alpha+(1-q_2)\beta)$.
Note that $\tilde\alpha-\tilde\beta=(q_1-q_2)^2(\alpha-\beta)$.
Then, if $q_1\not= q_2$ and $\alpha\not=\beta$, then 
$\tilde\alpha\not=\tilde\beta$.

\begin{thm}\label{conv_fi}
For $p(0,\cdot)=p^0\in L_+^1\cap L^\infty([0,1])$, the solution of 
Equation~(\ref{replicator_pde}) is unique, non-negative, and 
accumulates on the boundaries, i.e., $p^\infty\bydef\lim_{t\to\infty}p
=\pi_0[p^0]\delta_0+\pi_1[p^0]\delta_1$, where
$\pi_0[p^0]=1-\pi_1[p^0]$ and the {\em fixation probability} 
of $E_{q_1}$ strategists is given by
\begin{eqnarray*}
\pi_1[p^0]&=&\frac{\int_0^1\left[\int_y^1p^0(x)\d x\right]
\exp\left(-y^2(q_1-q_2)^2\frac{\alpha-\beta}{2}-
y(q_1-q_2)(q_2\alpha+(1-q_2)\beta)\right)\d y}
{\int_0^1\exp\left(-y^2(q_1-q_2)^2\frac{\alpha-\beta}{2}-
y(q_1-q_2)(q_2\alpha+(1-q_2)\beta)\right)\d y}\\
&=&
\frac{\int_0^1\int_0^xp^0(x)\exp\left(
-y^2(q_1-q_2)^2\frac{\alpha-\beta}{2}-y(q_1-q_2)(q_2\alpha+(1-q_2)\beta)\right)
\d y\ d x}
{\int_0^1\exp\left(
-y^2(q_1-q_2)^2\frac{\alpha-\beta}{2}-y(q_1-q_2)(q_2\alpha+(1-q_2)\beta)\right)
\d y}\ .
\end{eqnarray*}
\end{thm}

\begin{proof}
It is enough to prove for $q_1=1$ and $q_2=0$ (i.e.,
for Equation~(\ref{replicator_pde_pure})) and then
change the result from $(\alpha,\beta)$ to
$(\tilde\alpha,\tilde\beta)$.
Existence, non-negativeness and convergence to the
  boundaries follows  from Theorem~\ref{pde_prop}.

To obtain values $\pi_i[p^0]$, $i=1,2$,  
we multiply Equation~(\ref{replicator_pde_pure}) by $\psi(x)$
and integrate from 0 to 1. On assuming that $p$ is such that
  integration by parts can be performed and that no boundary terms
  arise, we obtain that
\[
\partial_t\int_0^1 p(t,x)\psi(x)\d x=
\int_0^1x(1-x)p(t,x)\left(\psi''(x)+(x(\alpha-\beta)+\beta)\psi'(x)\right)
\d x\ .
\]
Conservation laws are obtained solving 
$\psi''(x)+(x(\alpha-\beta)+\beta)\psi'(x)=0$. Solutions are
given by $\psi=\text{cte}$ (conservation of probability)
and 
\[
\psi(x)=c^{-1}\int_0^x\exp\left(-y^2\frac{\alpha-\beta}{2}-y\beta\right)\d
 y,\quad
\quad c=\int_0^1\exp\left(-y^2\frac{\alpha-\beta}{2}-y\beta\right)\d
 y\ .
\]

\begin{rmk}
Notice that $\psi(x)$ is the continuous counterpart to the discrete
fixation probabilities.
\end{rmk}

Using that 
\[
\int_0^1 p^0(x)\psi(x)\d x=
\int_0^1\left(\pi_0[p^0]\delta_0+\pi_1[p^0]\delta_1\right)\psi(x)\d x 
\]
we get
\begin{eqnarray*}
\pi_1[p^0]&=&\frac{\int_0^1\left[\int_y^1p^0(x)\d x\right]
\exp\left(-y^2\frac{\alpha-\beta}{2}-y\beta\right)\d y}
{\int_0^1\exp\left(-y^2\frac{\alpha-\beta}{2}-y\beta\right)\d y}\ ,\\
&=&\frac{\int_0^1\int_0^xp^0(x)\exp
\left(-y^2\frac{\alpha-\beta}{2}-y\beta\right)\d y\d x}
{\int_0^1\exp\left(-y^2\frac{\alpha-\beta}{2}-y\beta\right)\d y}\ .
\end{eqnarray*}
Finally, we change from $\alpha$, $\beta$ to $\tilde\alpha$, $\tilde\beta$.
\end{proof}

\begin{cor}
If $p^0=\delta_{x^0}$, then 
\[
\pi_1[\delta_{x^0}]=\frac{\int_0^{x_0}\exp
\left(-y^2(q_1-q_2)^2\frac{\alpha-\beta}{2}-
y(q_1-q_2)(q_2\alpha+(1-q_2)\beta)\right)\d y}
{\int_0^1\exp\left(-y^2(q_1-q_2)^2\frac{\alpha-\beta}{2}-
y(q_1-q_2)(q_2\alpha+(1-q_2)\beta)\right)\d y}\ .
\]
\end{cor}

\begin{deff}
We say that $E_{q_2}$ dominates $E_{q_1}$ ($E_{q_2}\succ E_{q_1}$)
if, for any initial condition $p^0\in L^1_+\cap L^\infty([0,1])$, 
the probability of fixation for the strategy $E_{q_1}$
is smaller than the one for 
the neutral case (the case $q_1=q_2$), i.e.,
\[
\pi_1[p^0]<\pi_1^\mathrm{N}[p^0]\bydef
\int_0^1 xp^0(x)\d x\ .
\]
We also say that $E_{q_2}$ $\delta$-dominates $E_{q_1}$ if
the above formula is valid for all $p^0=\delta_{x^0}$, $x^0\in(0,1)$,
i.e.,
\begin{equation}\label{weakdomination}
\pi_1[\delta_{x^0}]=\frac{\int_0^xF_{(q_1,q_2)}(y)\d y}
{\int_0^1F_{(q_1,q_2)}(y)\d y}< x^0\, \ \ \forall x^0\in(0,1)\ ,
\end{equation}
where we defined the auxiliary function
\begin{equation}\label{F_deff}
F_{(q_1,q_2)}(y)\bydef\exp\left(-y^2(q_1-q_2)^2\frac{\alpha-\beta}{2}
-y(q_1-q_2)(q_2\alpha+(1-q_2)\beta)\right)\ .
\end{equation}
\end{deff}
The following lemma shows that the two definitions above 
are in fact equivalent:

\begin{lem}
$E_{q_2}$ $\delta$-dominates $E_{q_1}$ if and only if $E_{q_2}\succ E_{q_1}$.
\end{lem}

\begin{proof}
We only need to prove the {\em only if} case.
Let us consider any initial condition given by 
$p^0\in L^1_+\cap L^\infty([0,1])$. Then
\begin{eqnarray*}
\pi_1[p^0]&=&\frac{\int_0^1\int_y^1 p^0(x)F_{(q_1,q_2)}(y)\d x\d y}
{\int_0^1F_{(q_1,q_2)}(y)\d y}\\
&=&
\frac{\int_0^1\int_y^1\int_0^1 p^0(z)\delta(z-x)F_{(q_1,q_2)}(y)\d z\d x\d y}
{\int_0^1F_{(q_1,q_2)}(y)\d y}\\
&=&
\frac{\int_0^1\int_0^1\int_y^x p^0(z)\delta(z-x)F_{(q_1,q_2)}(y)\d y\d z\d x}
{\int_0^1F_{(q_1,q_2)}(y)\d y}\\
&=&
\int_0^1\int_0^1p^0(z)\delta(z-x)\frac{\int_0^xF_{(q_1,q_2)}(y)\d y}
{\int_0^1F_{(q_1,q_2)}(y)\d y}\d z\d x\ .
\end{eqnarray*}
Now, we use Equation~(\ref{weakdomination}) and conclude
that 
\[
\pi_1[p^0]<\int_0^1\int_0^1 p^0(z)\delta(z-x)x\d z\d x=\int_0^1p^0(x)x\d x
=\pi_1^\mathrm{N}[p^0]\ .
\]
\end{proof}

In view of this lemma,
from now on, we consider only initial conditions of 
$\delta$-type, i.e., $p^0=\delta_{x^0}$.
In order to prove dominance relations, we prove first the
following:

\begin{lem}\label{increasing}
If $F_{(q_1,q_2)}$ is increasing in the interval $[0,1]$,
then $E_{q_2}\succ E_{q_1}$.
\end{lem}

\begin{proof}
For $p^0=\delta_{x^0}$, 
Equation~(\ref{weakdomination}) can be
re-written as
\[
\frac{1}{x^0}\int_0^{x^0}F_{(q_1,q_2)}(y)\d y<
\int_0^1F_{(q_1,q_2)}(y)\d y\ ,\ \ \forall x^0\in(0,1)\ .
\]
This equation can be interpreted as saying that the
average of the function $F_{(q_1,q_2)}$ in any interval
$[0,x^0]$, $x^0\in(0,1)$ is less than the average
in the interval $[0,1]$, which is true whenever 
the function is increasing.
\end{proof}

Finally we prove the full relations of dominance
for a $2\times 2$ game.
\begin{thm}
Let $E_{q_1}$ and $E_{q_2}$, $q_1,q_2\in[0,1]$, be two strategists in a 
$2\times 2$ game, and let $q^*=\beta/(\beta-\alpha)$. 
Then the relation of dominance is given by
Table~\ref{tabela2}.
\end{thm}

\begin{table}
\begin{center}
\begin{tabular}{c|c|c}
 &  & $E_{q_2}\succ E_{q_1}$ if and only if \\
\hline
$\alpha>\beta>0$&$q^*<0$&$q_2>q_1$\\
$\alpha>0>\beta$&$q^*\in(0,1)$&$q_2<q_1\le q^*$ or $q_2>q_1\ge q^*$\\
$0>\alpha>\beta$&$q^*>1$&$q_2>q_1$\\
$0>\beta>\alpha$&$q^*<0$&$q_2<q_1$\\
$\beta>0>\alpha$&$q^*\in(0,1)$&$q_1<q_2\le q^*$ or $q_1>q_2\ge q^*$\\ 
$\beta>\alpha>0$&$q^*>1$&$q_2>q_1$.
\end{tabular}
\caption{Dominance relations for the
non-degenerated ($\alpha\not=\beta\not=0\not=\alpha$) thermodynamical
limit of the frequency-independent Moran process, given by
Equation~\protect(\ref{replicator_pde}), with
$q^*=\beta/(\beta-\alpha)$.\label{tabela2}}
\end{center}
\end{table}

\begin{proof}
The proof consists in a long and tedious calculation
proving that, for each range in Table~\ref{tabela2},
the function $F_{(q_1,q_2)}$ is increasing. Then
we use Lemma~\ref{increasing}.
\end{proof}

\begin{figure}
\begin{center}
\epsfig{file=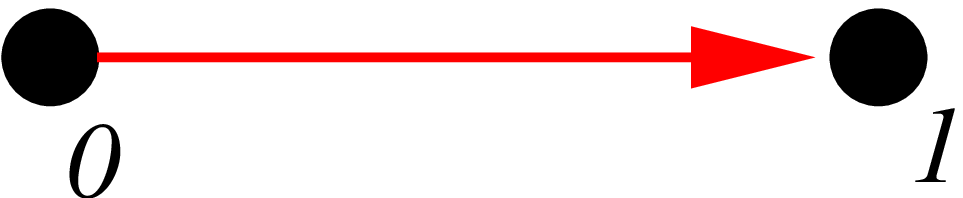,width=0.4\textwidth}
\hspace{.01\textwidth}
\epsfig{file=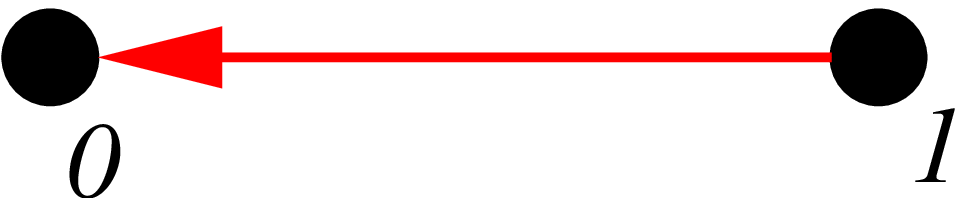,width=0.4\textwidth}
\hspace{.01\textwidth}
\epsfig{file=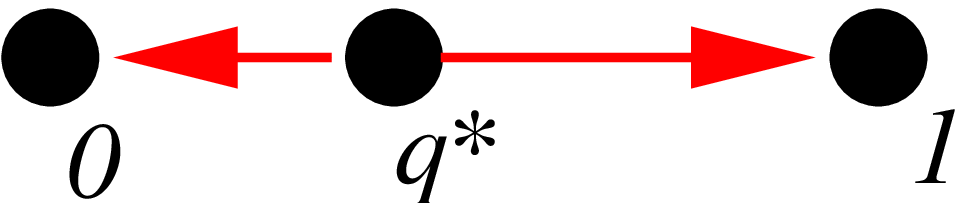,width=0.4\textwidth}
\hspace{.01\textwidth}
\epsfig{file=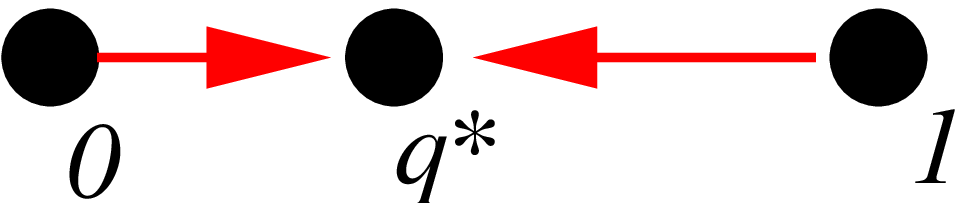,width=0.4\textwidth}
\vspace{.1cm}
\end{center}
\caption{Relation of dominance between $E_{q_1}$-- and 
$E_{q_2}$--strategists for given parameters. 
Here, $q^*=\beta/(\beta-\alpha)$.
The first 2 figures (above) show dominance from pure strategies,
the third one ($\alpha>0>\beta$, below, left)
shows that the pure strategies dominates their
neighbors (and everybody dominates $E_{q^*}$) and the last
one ($\beta>0>\alpha$, below, right)
shows that $E_{q^*}$ dominates any strategy.
The arrow points from the dominated to the dominant.
\label{dominance_fig}}
\end{figure}

The following corollary shows that the strategy $E_{q^*}$
is the best possible strategy if $\beta>0>\alpha$. 
\begin{cor}
If $\beta>0>\alpha$,
then $E_{q^*}\succ E_q$, $\forall q\not=q^*\bydef\beta/(\beta-\alpha)\in(0,1)$.
\end{cor}

\begin{proof}
For $q_2=q^*$, $F_{(q,q^*)}$ simplifies for 
\[
F_{(q,q^*)}(y)=\exp\left(-y^2(q-q^*)^2\frac{\alpha-\beta}{2}\right)\ .
\]
For $\alpha-\beta<0$, this is an increasing function of $y$ and this
proves the corollary.
\end{proof}

In order to finish the full picture of dominance, we need also
the following:
\begin{lem}\label{asymmetric_lem}
If $E_{q_2}\succ E_{q_1}$, then $E_{q_1}\not\succ E_{q_2}$.
\end{lem}

\begin{proof}
First, we see, from Equation~(\ref{F_deff}), that
\[
F_{(q_1,q_2)}(y)=F_{(q_2,q_1)}(1-y)F_{(q_1,q_2)}(1)\ .
\]
Then, we write
\[
\frac{1}{x^0}\int_0^{x^0}F_{(q_1,q_2)}(y)\d y=
\frac{1}{x^0}\left[\int_0^1F_{(q_2,q_1)}(y)\d y-\int_0^{1-x^0}F_{(q_2,q_1)}(y)\d y
\right]F_{(q_1,q_2)}(1)\ .
\]
Furthermore,
\[
\int_0^1F_{(q_1,q_2)}(y)\d y=\int_0^1F_{(q_2,q_1)}(y)\d y F_{(q_1,q_2)}(1)\ .
\]
Using the fact that $E_{q_2}\succ E_{q_1}$, we find
\[
\frac{1}{x^0}\left[\int_0^1F{(q_2,q_1)}(y)\d y-\int_0^{1-x^0}F_{(q_2,q_1)}(y)
\right]<\int_0^1F_{(q_2,q_1)}(y)\d y\ .
\]
We re-arrange the terms and conclude that
\[
\frac{1}{1-x^0}\int_0^{1-x^0}F_{(q_2,q_1)}(y)\d y>\int_0^1F_{(q_2,q_1)}(y)\d y
\ .
\]
\end{proof}

Table~\ref{tabela2}, together with Lemma~\ref{asymmetric_lem} 
can be summarized in Figure~\ref{dominance_fig}.

It is important to note also that in some 
references (see, e.g.,~\cite{TaylorSasaki})
it is said that selection favors strategy II replacing
strategy I (
in this case, we say that strategy II
{\em weakly dominates} strategy I), in a finite population
of size $N$, if a single type II mutant has fixation 
probability larger than $1/N$, the neutral probability.
Unfortunately, no sound generalization of
this concept can have a graph similar to the one presented
in Fig.~\ref{dominance_fig}, as it is possible that $E_{q_2}$ weakly
dominates strategy $E_{q_1}$ and
vice-versa. See~\cite{TaylorSasaki} for details.

The concept explained above clearly extend the concept
of ESS for the PDE case. As the PDE case works as an 
approximation for large $N$ of the discrete case,
it is easy to see that we can extend the ESS definition
also to the more realistic discrete case. 

Different to the most known ODE (see also the next
section) case for the definition
of ESS, here we cannot guarantee that the probability
distribution
will, in the long range (when adequately parametrized) 
accumulate in the ESS (when it is in the interior
of the interval $[0,1]$), but we can see that an individual
that plays strategy I and II with frequencies given by
the game's ESS is optimized to win any contest (with the
same parameters).

If the strategists involved in the game play with
frequencies different from the ESS (for example, the
pure strategies) the ODE prediction is that a stable
mixture will evolve. This is impossible in the
discrete case (as, in the long range, all individuals
will descend of a single one in time $t=0$, which will be
of one of the given types) and also in the PDE model (as shown by
Theorem~\ref{conv_fi}).

More generally, we say

\begin{thm}
Let $p_{N,\Delta t}(x,t)$ be the solution of the finite population
dynamics (of population $N$, time step $\Delta t=1/N^2$), with initial
conditions given by $p^0_N(x)=p^0(x)$, $x=0,1/N,2/N,\cdots,1$, for
$p^0 \in L^1_+([0,1])$. Assume also that $(A-1,B-1,C-1,D-1)=
1/N(a,b,c,d)+\O(1/N^2)$.
Let $p(t,x)$ be the solution of the continuous model
with initial condition given by $p^0(x)$.
If we write $p_i^n$ for the $i$-th component of $p_{N,\Delta t}(x,t)$
in the $n$-th iteration, we have, for any $t^{*}>0$, that
\[
\lim_{N\to\infty}  p^{tN^2}_{xN}
=p(t,x),\quad x\in[0,1],\quad t\in[0,t^{*}].
\]
\end{thm}

\begin{proof}
First, we consider the matrix $\widetilde{\mat}$ obtained from $\mat$ by
deleting the first and last rows and columns.
Then, we observe that the derivation of the thermodynamical limit
shows that the discrete iteration given by $\widetilde{\mat}$ is consistent --- in
the approximation sense~\cite{richtmyer:morton} --- with  
Equation (\ref{replicator_pde_pure}), without any boundary conditions,
provided that we set $A=1+a/N$, and similarly  for $B$, $C$ and
$D$. From the results of Appendix A, we know that the discrete
iteration is stable, since $\sigma(\widetilde{\mat})\subset(-1,1)$. From
Appendix B, we see that the continuous problem without boundary
conditions is well posed in the $D_s$ spaces defined there. In this
case, we can then invoke  the Lax-Ricthmyer equivalence theorem
\cite{richtmyer:morton} to guarantee that the discrete model converges to 
the continuum one, in the limit $\Delta t,\Delta x\to0$, with $\Delta
t=(\Delta x)^ 2$. More precisely, the iteration defined by
$\widetilde{\mat}$ converges to $q(t,x)$, the smooth part of $p(t,x)$;
cf. appendix~\ref{app_b}

Now returning to the iteration defined by $\mat$. In order to finish
the proof, we only need to show
that $P(t,0)$ and $P(t,1)$ converges weakly to the appropriate Dirac
masses. We shall do the computation for $x=0$, the case $x=1$ being
similar.

For $x=0$ the iteration defined by $\mat$ reads
\[
P(t+\Delta t,0)=P(t,0)+\frac{1}{N}P\left(t,\frac{1}{N}\right)
\]
Thus, letting $t=0$ and solving the recursion, we have that
\[
P(m\Delta t,0)=P(0,0)+\frac{1}{N}\sum_{j=1}^{m-1}P\left(j\Delta t,\frac{1}{N}\right).
\]
Since $\mathbf{e}_1\in\mathbb{R}^N$ converges weakly to
$\delta_0$ as $N\to\infty$---by considering test functions with
support contained in $(1/N,1/N)$---we need only to show that it 
has the correct mass at each time $t$. For this, notice that 
\[
\mathcal{P}\left(j\Delta t,\frac{1}{N}\right)=\mathcal{P}\left(j\Delta
t,0\right)+\frac{1}{N}\partial_x\mathcal{P}\left(j\Delta t,0\right).
\]
Since $p(t,x)=N\mathcal{P}(t,x)$, we find that, in a weak sense,
\[
\lim_{N\to\infty}p^{tN^2}_0\to\int_0^tq(s,0)\d s+ P(0,0).
\]
\end{proof}

\section{The diffusionless case and the replicator
dynamics}
\label{diffusionless}

We shall see in this Section that the ODE Replicator dynamics is
  equivalent to the diffusionless version of
  Equation~(\ref{replicator_pde}). This will have important
  consequences that we shall discuss later on. Notice also,
  cf. Remark~\ref{scalling_rmk}, that this is the correct limiting
  equation, if the payoffs decay slowly to one as $N\to\infty$.

Thus, we consider
\begin{equation*}
\partial_t p=-\partial_x\left(x(1-x)(x(a-c)+(1-x)(b-d))p\right)\ .
\end{equation*}
A weak solution of this equation is given by 
$p(t,x)=\delta_{X(t)}$,  if $X(t)$ solves
\begin{equation}
\label{replicator_ode}
\dot X=X(1-X)(X(a-c)+(1-X)(b-d))\ ,
\end{equation}
which is the simplest replicator equation~\cite{SigmundHofbauer} 
for the two-person game with pay-off
matrix given by
\begin{equation}\left(
\begin{matrix}
a&b\\
c&d
\end{matrix}
\right)\ .
\end{equation}

The stationary points of Equation~(\ref{replicator_ode}) are
given by $0$, $1$ and $X^*\bydef\beta/(\beta-\alpha)$. The most
interesting scenario occurs when $\alpha<\beta$ and $X^*\in(0,1)$
(i.e, $\beta>0>\alpha$):
in this case the only stable equilibrium is the non-trivial 
$X=X^*$. For the full analyze, see Table~\ref{tabela1}.
Compare also with the description of dominance in the 
previous section.

\begin{table}
\begin{center}
\begin{tabular}{c|c||c|c}
 & & stable & unstable \\
\hline
$\alpha>\beta>0$&$X^*<0$&$X=1$&$X=0$\\
$\alpha>0>\beta$&$X^*\in(0,1)$&$X=0$ and $1$&$X=X^*$\\
$0>\alpha>\beta$&$X^*>1$&$X=0$&$X=1$\\
$0>\beta>\alpha$&$X^*<0$&$X=0$&$X=1$\\
$\beta>0>\alpha$&$X^*\in(0,1)$&$X=X^*$&$X=0$ and $1$\\
$\beta>\alpha>0$&$X^*>1$&$X=1$&$X=0$
\end{tabular}
\caption{Stable and unstable equilibria in
the range $[0,1]$  for the 
non-degenerated ($\alpha\not=\beta\not=0\not=\alpha$) replicator
dynamics~\protect(\ref{replicator_ode})}\label{tabela1}
\end{center}
\end{table}

\begin{figure}
\label{flux_fig}
\begin{center}
\epsfig{file=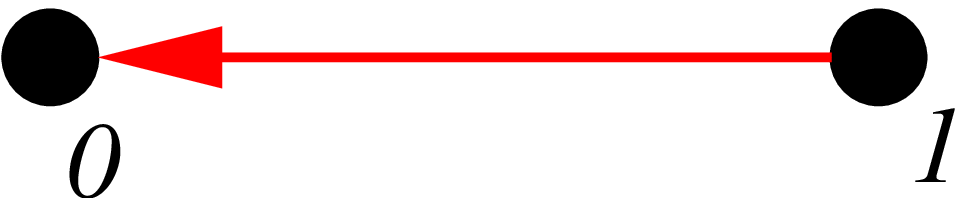,width=0.4\textwidth}
\hspace{.01\textwidth}
\epsfig{file=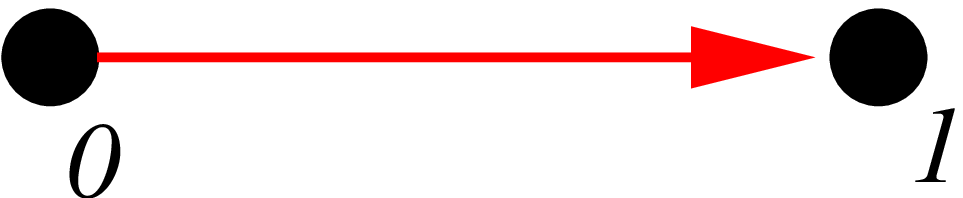,width=0.4\textwidth}
\hspace{.01\textwidth}
\epsfig{file=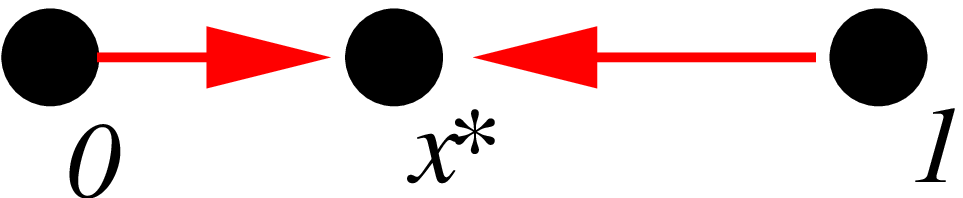,width=0.4\textwidth}
\hspace{.01\textwidth}
\epsfig{file=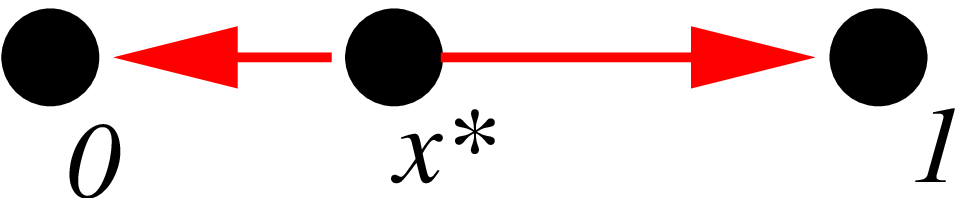,width=0.4\textwidth}
\vspace{.1cm}
\end{center}
\caption{Flux of the replicator equation for 
pure strategy dominated games (above), mixed strategy
dominated (below, left) and bistable games (below, right). 
Here, $X^*=\beta/(\beta-\alpha)$.
Compare with Figure~\protect\ref{dominance_fig}.}
\end{figure}

Our definition of dominance seems more general than
many definitions that appear in the 
literature~\cite{Neill1,NowakSazaki1,Scaffer1,TaylorSasaki}. 
Furthermore, the use of the thermodynamical limit in the analysis 
make it much more simple to work. In particular, consider
a game between $E_{q_1}$- and $E_{q_2}$-strategists 
and a given replicator dynamics
such that any non-trivial initial conditional 
converges in $t\to\infty$ to one of the two
trivial equilibria, say, $X=0$. 
The replicator dynamics is given by
\begin{equation}\label{full}
\dot X=X(1-X)
(X(q_1-q_2)^2(\alpha-\beta)+(q_1-q_2)(q_2\alpha+(1-q_2)\beta))
\ .
\end{equation}
If, for any initial condition $X(0)\in(0,1)$, 
$\lim_{t\to\infty}X(t)=0$, then, $\dot X(t)<0$, $\forall X\in(0,1)$,
$\forall t\in\R^+$, 
i.e., $(X(q_1-q_2)^2(\alpha-\beta)+(q_1-q_2)(q_2\alpha+(1-q_2)\beta))
<0$. This implies that $F'_{(q_1,q_2)}(y)>0$, $\forall y\in(0,1)$,
where $F_{(q_1,q_2)}$ is defined by Equation~(\ref{F_deff}). 
In particular $F_{(q_1,q_2)}$ is increasing and then
from Lemma~\ref{increasing},
we have $E_{q_2}\succ E_{q_1}$. If, on the other hand,
$\lim_{t\to\infty}X(t)=1$, by a similar argument, we
have that $E_{q_1}\succ E_{q_2}$.

These picture is completed after looking to
Figures~\ref{dominance_fig} and~\ref{flux_fig} and noting that
the flow of the replicator dynamics always goes
from the less-dominant strategy to the more dominant
one, if we consider an equivalence (at the replicator
dynamics level) between mixed populations of 
pure strategist and populations of mixed strategists.

In reference~\cite{TraulsenClaussenHauert} a thermodynamical
limit of a frequency-dependent Moran process was also designed, 
but
the pay-off were not re-scaled when $N\to\infty$ and the
Fokker-Planck equation obtained was claimed to be valid for 
{\em large, but finite} $N$, and not in the thermodynamical
limit.

\section{The Frequency Independent Moran Process}
\label{moran}

In order to consider frequency-independent fitness,
we impose a pay-off matrix such that the 
gain of a player is independent of others player's
strategies, that is, $A=B$ and $C=D$. In particular,
we impose $C/A=D/B=r$. The number $r$ is know as
the relative fitness. Most results here are simple
corollaries of results from the previous section. 
We state them only for completeness.

\begin{cor}\label{cor:fixation}
The fixation probabilities $F_n$ of type $\mathbb{A}$ individuals
for an initial condition of $n$ mutants
in the frequency independent Moran process with
relative fitness $r$ are given
by
\begin{eqnarray}
\label{r_not_1}
&&F_n=\frac{1-r^n}{1-r^{1-N}}+\frac{k}{N}\frac{r^n-r^{1-n}}{1-r^{1-N}}\ ,\\
\label{r_1}
&&F_n=\frac{n}{N}\ ,\ \ r=1\ .
\end{eqnarray}
\end{cor}

\begin{proof}
When the relative fitness is constant, 
i.e. $\rho_N(n)= 1/r$, (\ref{gen:fixsoln}) becomes
\begin{align}
F_n &=
G_1\sum_{k=1}^n\frac{1}{\rho^{k-1}}\left(1+\frac{(\rho-1)(k-1)}{N-1}\right),
\nonumber \\
G_1&=\left[\sum_{k=1}^N\frac{1}{\rho^{k-1}}
\left(1+\frac{(\rho-1)(k-1)}{N-1}\right)
\right]^{-1}
\label{crgen:fixsoln}
\end{align}

We 
sum the series and prove the corollary.
If $r=1$, it is straightforward to see that
$F_n=n/N$.  

\end{proof}

\begin{rmk}\label{rmk:bddb}
In the case of birth/death process we have instead:
\[
F_n=G_1\sum_{k=1}^n\frac{1}{\rho^{k-1}}=G_1\sum_{k=1}^nr^{k-1}=G_1\frac{1-r^k}{1-r}
=\frac{1-r^k}{1-r^N}\ ,
\]
where we used that
\[
G_1=\frac{1-r}{1-r^N}.
\]
\end{rmk}

Note that the coefficients obtained in Corollary~\ref{cor:fixation}
are different from the
one obtained in~\cite{Komarova1}, which are the same
as in Remark~\ref{rmk:bddb}. The difference is the result
of differences between a death/birth and birth/death processes.
Anyhow, the formulas
are equivalent for large $N$.

We also define $\gamma\bydef\alpha=\beta$ and then
Equation~(\ref{replicator_pde}) is
\begin{equation}\label{our_equation}
\partial_tp=\partial_x^2\left(x(1-x)p\right)-
\gamma\partial_x\left(x(1-x)p\right)\ .
\end{equation}

As a simple consequence of Theorem~\ref{conv_fi} for
Equation~(\ref{our_equation}), we have
\begin{cor}
Let $p$ be a solution of Equation~(\ref{our_equation}) with initial
conditions $p^0\in L^1_+\cap L^\infty([0,1])$. Then, 
in a weak sense, $p^\infty\bydef\lim_{t\to\infty}p(\cdot,t)=
\pi_0[p^0]\delta_0+\pi_1[p^0]\delta_1$. 
Furthermore, we have $\pi_0[p^0]=1-\pi_1[p^0]$
and
\[
\pi_1[p^0]=\frac{1-\int_0^1\e^{-\gamma x}p^0(x)\d x}{1-e^{-\gamma}}\ .
\]
\end{cor}
If we start with $p^0=\delta_{x_0}$, then 
\begin{equation}
\label{fix_single}
f_\gamma(x_0)\bydef
\pi_1[\delta_{x_0}]=\frac{1-\e^{-\gamma x_0}}{1-\e^{-\gamma}}\ ,
\end{equation}
and $\lim_{\gamma\to 0}\pi_1[\delta_{x_0}]=x_0$. 
This is true because the neutral case corresponds to
$\gamma=0$.
Note that
$f_\gamma(0)=0$, $f_\gamma(1)=1$, $\forall\gamma$ and
that $f_\gamma(x_0)\ge x_0$ if and only if $\gamma\ge 0$. So, in the
language of previous sections, $\mathbb{A}\succ\mathbb{B}\iff \gamma>0$.

It is important to compare the probability of
fixation in the continuous limit, Equation~(\ref{fix_single}),
and the result obtained for finite population, 
Equation~(\ref{r_not_1}). To understand the idea we should
consider that, in the finite case, we have
initially a fixed proportion $\kappa\in(0,1)$ of mutants, 
such that the probability of fixation is given by
\[
\frac{1-r^{\kappa N}}{1-r^{N-1}}-
\kappa\frac{r^{\kappa N}-r^{\kappa N-1}}{1-r^{N-1}}
\approx
\frac{1-r^{\kappa N}}{1-r^N}\ ,
\]
when $N$ is large and $r$ close to 1. To be more precise,
if $r(N)=1+\gamma/N$,
\[
\lim_{N\to\infty}
\frac{\frac{1-r^{\kappa N}}{1-r^{N-1}}-
\kappa\frac{r^{\kappa N}-r^{\kappa N-1}}{1-r^{N-1}}}
{\frac{1-r^{\kappa N}}{1-r^N}}=1\ .
\]
In order to compare that formula with~(\ref{fix_single}) for
large $N$, we need only to impose $k=x_0$ (the initial fraction
of mutants) and then $\e^{-\gamma}\approx r^{-N}$, i.e., 
$\gamma\approx N(r-1)$, for $r-1\ll 1$ (valid for large $N$),
in agreement with $\gamma=\lim_{N\to\infty}N(r-1)$ (compare 
with~(\ref{Ntoinfty2})).

We cannot avoid the comparison of our result with the
classical results by Kimura~\cite{Kimura}. Following 
this reference, let $u(t,y)$ be the probability
that a mutant allele, initially with frequency $y$
and relative fitness $s$ be fixed after a time $t$
in a randomly mating diploid population of size $N_0$.
Then
\begin{equation}\label{Kimura_eq}
\partial_t u=\frac{y(1-y)}{4N_0}\partial_y^2u+
sy(1-y)\partial_y u\ .
\end{equation}

This equation and Equation~(\ref{our_equation}) are
 associated backward/forward Kolmogorov equations 
with suitable rescalings~\cite{Gardiner}. 
Then, for example, Equation~(\ref{fix_single})
is the same found in~\cite{Kimura}, where $\gamma=4Ns$
for $s=r(N)-1$ is the selective advantage.
Furthermore,
the fact that $f_0(x_0)=x_0$ reproduces the idea that
a neutral mutant ($\gamma=0$) is fixed with probability 
equal to its initial frequency. 

Following, again, reference~\cite{Gardiner}, if $u(t,y)$ solves 
Equation~(\ref{Kimura_eq}), then $u(t,x)=u_\mathrm{S}(x)p(x,t)$
where 
\[
u_\mathrm{S}(x)=\frac{1-\e^{-N_0sx}}{1-\e^{-4N_0s}}
\]
is the stationary solution of Equation~(\ref{Kimura_eq}) and
$p(x,t)$ solves~(\ref{our_equation}) (with appropriate rescalings
and normalizations). This shows the
equivalence of this deduction and Kimura's one.

\section{The drift limit}
\label{drift}

The ``drift limit'' means that the time-step is
re-scaled according to $\Delta t=1/N$. In this case, we do not need to
consider the weak selection limit, i.e., pay-offs (and fitness)
are considered time-step independent.
This problem is mathematically well posed, but, as explained
below, it seems not to be an interesting limit from the modeling
point of view. We state it only for completeness.

First, we see what happens for the drift limit of the
frequency dependent Moran process, i.e., 
Equation~(\ref{drift-independent}).

\begin{thm}
Let $\bar p$ be the solution of Equation~(\ref{drift-independent})
with initial conditions given by $\bar p^0\in L^1_+\cap L^\infty([0,1])$.
Then, 
$\bar p^\infty=\bar\pi_0\delta_0+\bar\pi_*\delta_{x^*}+\bar\pi_1\delta_1$,
where $\bar\pi_0+\bar\pi_*+\bar\pi_1=1$ 
and $x^*=-(B-D)/(A-B-C+D)$. Furthermore,
if $A-C<0$, then $\bar\pi_0=0$; if $B-D>0$ then $\bar\pi_1=0$; and
if $(AD-BC)/((A-C)(B-D))<0$ then $\bar\pi_*=0$.
If $x^*\not\in[0,1]$, $\bar\pi_*=0$.
\end{thm}

\begin{proof}
We multiply Equation~(\ref{drift-independent}) by
\[
\psi(x)=(1-x)^{A/(A-C)}x^{-D/(B-D)}(x(A-B-C+D)+B-D)^{(DA-BC)/((A-C)(B-D))}
\]
and integrate from 0 to 1. Then
\begin{eqnarray*}
\partial_t\int_0^1\psi(x)\bar p(x,t)\d x&=&
-\int_0^1
\frac{x(1-x)(x(A-B-C+D)+B-D)}{x^2(A-B-C+D)+x(B+C-2D)+D}\psi'(x)\bar p(x,t)\d x 
\\
&=&-\int_0^1\psi(x)\bar p(x,t)\d x\ .
\end{eqnarray*}
From Gronwall's inequality, we find that $\bar p^\infty$ is
supported at the zeros of $\psi(x)$. 
\end{proof}

Suppose that we have a game where the strategy I dominates
(e.g., the Prisoner's dilemma, where strategy I means
``defect''),
i.e., $A>C$ and $B>D$. If $AD-BC>0$, $\bar\pi_*=0$,
and if $AD-BC<0$, then $x^*>1$, and this implies $\bar\pi_*=0$.
Eventually, the full population will play strategy I.

For $A<C$ and $B<D$, the full population will play strategy II.
 
For the Hawk-and-Dove game we have $A-C<0$ and $B-D>0$.
This implies that $(AD-BC)/((A-C)(B-D))>0$ and then 
$\bar p^\infty=\delta_{x^*}$, where $x^*\in(0,1)$.

Finally, for coordination games, 
$A-C>0$ and $B-D<0$, then $(AD-BC)/((A-C)(B-D))<0$
and $\bar p^\infty=\bar\pi_0\delta_0+\bar\pi_1\delta_1$. To
obtain the values $\bar\pi_i$, $i=0,1$, note that
$x^*\in(0,1)$ and
\[
\partial_t\int_0^{x^*}\bar p\d x=0\ ,\ \ \partial_t\int_{x^*}^1\bar p\d x =0.  
\]
This implies that
\begin{eqnarray*}
\bar\pi_0&=&\int_0^{x^*}\bar p^\infty\d x=\int_0^{x^*}\bar p^0\d x\ ,\\
\bar\pi_1&=&\int_{x^*}^1\bar p^\infty\d x=\int_{x^*}^1\bar p^0\d x\ .
\end{eqnarray*}
In a pictorial way, all the mass to the right of $x^*$ will move toward 
the point $x=1$, while the mass on the left will move toward 0.
If the initial condition is of delta-type, i.e., $p^0=\delta_{x^0}$
then the final condition is fully determined, $\bar p^\infty=\delta_0$
($\bar p^\infty=\delta_1$) if $x_0<x^*$ ($x_0>x^*$, respectively).

Now, we consider the frequency independent case, i.e., we
impose $A=B=1$ and $C=D=r$ at Equation~(\ref{drift-independent}). 

\begin{cor}\label{r_indp_delta_t}
Let $\bar p$ be the solution of 
\begin{equation}\label{Moran_falso}
\partial_t\bar p=
-(r-1)\partial_x\left[\frac{x(1-x)}{x(r-1)+1}\bar p\right]\ .
\end{equation}
with $\bar p^0\in L^1_+\cap L^\infty([0,1])$.
Then $\bar p^\infty=\delta_1$ for
$r>1$ and $\bar p^\infty=\delta_0$ for $r<0$.
\end{cor}

\begin{proof}
Note that $\psi(x)=(1-x)^{1/(1-r)}x^{-r/(1-r)}$. Then, its
zeros are at most 0 and 1. This implies $\bar\pi_*=0$. 
The values of $\bar\pi_0$ and $\bar\pi_1$ follow 
trivially.
\end{proof}

As a conclusion of this corollary, we note that the time-step
of order $1/N$ implies in no diffusion, i.e., no genetic
drift. So, the result of Equation~(\ref{Moran_falso}) is
deterministic, in the sense that an arbitrarily small fraction
of advantageous mutant will eventually take over the 
entire population, while disadvantageous mutants will
certainly be extinct (if the population is initially
mixed). In Equation~(\ref{our_equation}) nothing
similar happens.

\section{Final Remarks}
\label{final}

The procedure used here can be applied to different
evolution process. For example, consider the
{\em imitation dynamics} given by the following 
rules: from a population with size $N$ and
two possible types, we choose two individuals
$I_1$ and $I_2$.
If they are of the same type, nothing changes.
If $I_1$ is of type $\mathbb{A}$ and $I_2$ of type
$\mathbb{B}$, $I_1$ changes its type with
probability $\Psi(\phi_B-\phi_A)$ and the same 
if we swap $I_1$ and $I_2$, where $\phi_A$ and $\phi_B$
are the fitness for the types $\mathbb{A}$ and $\mathbb{B}$
respectively and $\Psi:\R\to[0,1]$ is a continuously
differentiable non decreasing function. 
Then, the
transition coefficients are given by
\begin{eqnarray*}
c_+(n,N)&=&\frac{N-n}{N}\frac{n}{N-1}\Psi(\phi_A-\phi_B)\ ,\\
c_-(n,N)&=&\frac{n}{N}\frac{N-n}{N-1}\Psi(\phi_B-\phi_A)\ ,\\
c_0(n,N)&=&1-c_+(n,N)-c_-(n,N)\ .
\end{eqnarray*}
We consider the functions of $x=n/N$ as defined 
in~(\ref{exp_c+})--(\ref{exp_c-}) and with 
assumptions~(\ref{Ntoinfty1})--(\ref{Ntoinfty2})
we get
\begin{eqnarray*}
\lim_{N\to\infty}N\left( c_1^{(1)}+c_0^{(1)}+c_-^{(1)}\right)&=&
6x^2\Psi'(0)(a-b-c+d)\\
&&+4x\Psi'(0)(-a+2b+c-2d)-2\Psi'(0)(b-d)\ ,\\
\lim_{N\to\infty}N\left(c_+^{(0)}-c_-^{(0)}\right)&=&-2x^3\Psi'(0)(a-b-c+d)\\
&&-2x^2\Psi'(0)(-a+2b+c-2d)+2x\Psi'(0)(b-d)\ ,\\
\lim_{N\to\infty}\left(c_+^{(2)}+c_0^{(2)}+c_-^{(2)}\right)&=&-4\Psi(0)\ ,\\
\lim_{N\to\infty}\left(c_+^{(1)}-c_-^{(1)}\right)&=&2\Psi(0)(2x-1)\ ,\\
\lim_{N\to\infty}\left(c_+^{(0)}+c_-^{(0)}\right)&=&2x(1-x)\Psi(0)\ .
\end{eqnarray*}
Gathering everything in Equation~(\ref{EqForp})
we find as the drift-diffusion limit of this process
\begin{equation}
\label{continuation}
\partial_tp=\Psi(0)\partial_x^2\left(x(1-x)p\right)
-2\Psi'(0)\partial_x\left(x(1-x)(x\alpha+(1-x)\beta)p\right)
\ ,
\end{equation}
with $\alpha=a-c$ and $\beta=b-d$.
From the assumptions, $\Psi(0),\Psi'(0)\ge 0$.
Relation of dominance for $E_{q_1}$- and $E_{q_2}$-strategists
are exactly the same as before, as can be easily computed from the fact
that the conservation laws associated to Equation~(\ref{continuation}) 
are $\psi(x)=1$ and 
\begin{equation*}
\psi(x)=\int_0^x\exp\left(\left(-\frac{y^2}{2}(q_1-q_2)^2(\alpha-\beta)
+y(q_1-q_2)(q_2\alpha+(1-q_2)\beta)\right)\frac{\Psi'(0)}{\Psi(0)}\right)
\d y\ .
\end{equation*}

The coefficients can be adjusted from the basic discrete
process. In particular, we can choose $\Psi$ such that
Equation~(\ref{continuation}) is drift-dominated
(if $\Psi(0)\ll \Psi'(0)$) or diffusion-dominated (if
$\Psi'(0)\ll \Psi(0)$). In a forthcoming paper, we will
completely study this equation and this two different
regimes. In particular,
we can define a family of functions $\Psi_\eps$, such that
$\lim_{\eps\to 0}\Psi_\eps(0)=0$, but $\lim_{\eps\to 0}\Psi_\eps'(0)>0$
and use singular-perturbation theory to understand the
diffusionless limit of the replicator-diffusion 
equation~(\ref{continuation}).
We can expect a behavior similar to the one
found in Section~\ref{diffusionless}. This means
that, for certain imitation dynamics and for intermediate times,
the evolution of the system, or more precisely, the 
``peak'' of the density distribution, can be modeled by
Equation~(\ref{replicator_ode}), as we can see in
Figures~\ref{mmmmmfig} and \ref{mmmmmmfig}.

\section*{Acknowledgments}

FACCC had his research supported by Project POCI/MAT/57546/2004.
MOS thanks Milton Lopes Filho for helpful discussions.

\appendix
\section{Proof of Proposition~\ref{Akto0}}
\label{ProofProp1}

The following properties of $f_N$ and $g_{N}$
  will be useful in the sequel:
\begin{enumerate}
\item $f_N(0)=f_N(N)=0$;
\item $g_N(0,r)=1-1/N$ and $g_N(N,r)=r-1/N$;
\end{enumerate}
Notice that the first column of $\mat$ is $e_1$ and last one is
$e_{N+1}$. Hence $1\in\sigma(\mat)$, and $e_1,e_{N+1}$ are associated
eigenvectors. Also, since $\mat$ is nonnegative tridiagonal, we must
have $\sigma(\mat)\subset\mathbb{R}$.

Since $\mat$ is column stochastic, and since all diagonal elements are
nonzero, an application of Gersgorin theorem
to $\mat^{\dagger}$ shows that $\sigma(\mat)\subset(-1,1]$.

We now show that the 1 is an eigenvalue of multiplicity two. 
First, observe that $\mat$ has the following block structure
\[
\begin{pmatrix}
1&*\hfill&\\
&\widetilde{\mat}&\\
&\hfill *&1\\
\end{pmatrix},
\]
where $\widetilde{\mat}$ is a $(N-1)\times (N-1)$ tridiagonal matrix, with
nonzero elements in the super and subdiagonal. Hence, $\widetilde{\mat}$ is 
irreducible. 

Let $\eta_i$ denote the sum of elements of the $i$-th column of
$\widetilde{\mat}$. Then we have
\[
\eta_i=1, i=2,\ldots,N-2\quad\text{and}\quad 0<\eta_1,\eta_{N-1}<1.
\]
Because of the irreducibility of $\widetilde{\mat}$, the strict inequality
for $\eta_1$ (or $\eta_{N-1}$) is sufficient to show that
$1\not\in\sigma(\widetilde{\mat})$ (cf. \cite{HJ}).

This result on the spectrum of $\widetilde{\mat}$, together with the block
structure of $\mat$ proves the claim.

We write
\[
\mat=P\Lambda P^{-1},
\]
where
\[
P=\begin{pmatrix}
1&***&0\\
\vdots&***&\vdots\\
0&***&1
\end{pmatrix}
\quad\text{and}\quad
\Lambda=\begin{pmatrix}
1&0\ldots0&0\\
0&J&0\\
0&0\ldots0&1\\
\end{pmatrix}
\]
We also notice that $P^{-1}$ has the same structure of $P$.

From the localization results on eigenvalues of $\mat$, we know that 
$\sigma(J)\subset(-1,1)$, and hence
\[
\lim_{k\to\infty}J^k=0.
\]

In this case, we have that:
\[
\lim_{k\to\infty}\Lambda^k=\begin{pmatrix}
1&0&\ldots&0\\
\vdots&\vdots&\vdots&\vdots\\
0&\ldots&0&1
\end{pmatrix},
\]
and the result follows.

\section{Proof of Theorem~\ref{pde_prop}}
\label{app_b}

First, if we let $\alpha=a-c$ and $\beta=b-d$ in
(\ref{replicator_pde}), we have
\begin{equation}
\partial_tp=\partial^2_x[x(1-x)p]-\partial_x[x(1-x)(\beta
  +(\alpha-\beta)x)p].
\label{wrk_pde}
\end{equation}

Further, let 
\[
\e^{\frac{1}{2}\left(\beta x+(\alpha-\beta)\frac{x^2}{2}\right)}w(t,x)=x(1-x)p(t,x).
\]
Equation (\ref{wrk_pde}) then becomes
\begin{equation}
\partial_tw=x(1-x)\left\{\partial^2_xw-\left[\frac{\alpha-\beta}{2}+\frac{1}{4}\left(\beta+(\alpha-\beta)x\right)^2\right]w\right\}.
\label{sawrk_pde}
\end{equation}

First, we observe that by writing $w_{\epsilon}=w+\epsilon t$ allows
us to prove a maximum principle for $C^2(\mathbb{R}^+\times(0,1))$ solutions to
(\ref{sawrk_pde}) in a standard way. In particular, since $w\geq0$ is in
the parabolic boundary, it is nonnegative everywhere.

Existence can be established by Fourier series theory. In what
follows, all the Banach spaces in this section are weighted with respect to 
\begin{equation}
\omega(x)=\frac{1}{x(1-x)}
\label{the:weight}
\end{equation}
Consider the associated equation
\begin{equation}
-\psi''+\left[\frac{\alpha-\beta}{2} +
 \frac{1}{4}\left(\beta+(\alpha-\beta)x\right)^2\right]\psi=\lambda
 \omega(x)\psi,\quad 
 \psi(0)=\psi(1)=0.
\label{our_eigenop}
\end{equation}
Since $w(x)\in L^1_{\mathrm{loc}}((0,1))$, standard Liouville theory
applies to (\ref{our_eigenop}). The relevant facts are collected in

\begin{lem}
Equation (\ref{our_eigenop}) defines a singular Sturm-Liouville problem
satisfying the following:
\begin{enumerate}
\item The extreme points are singular points of limit point,
  non-oscillatory type. The  Friedrich's extension of the operator on 
  the left hand side of (\ref{our_eigenop}) is a self-adjoint operator
  in $H^2([0,1])\cap H^0([0,1])$,  that is bounded from below.
\item The eigenvalues of (\ref{our_eigenop}) are real, purely
  discrete, bounded from  below, and accumulate only at infinity.
\item The associated eigenfunctions are an orthonormal basis of
  $L^2([0,1])$.
\item If $\{\lambda_j\}$ denotes the spectrum, we have
\[
\lim_{j\to\infty}\frac{\lambda_j}{j^2}=K\not=0.
\]
\end{enumerate}
\end{lem}

\begin{proof}

A straightforward Frobenius analysis near 0 and 1, shows that only one
of the linear independent solutions can square integrable with respect
to $\omega(x)$. Moreover, the Frobenius expansion are regular without complex
exponents. Hence, the extremes are of limit point, non-oscillatory
type. The other results are standard---see 
for instance \cite{CL}.
\end{proof}
An important property of (\ref{our_eigenop}) is given by
\begin{lem}
The operator defined by (\ref{our_eigenop}) is positive-definite.
\end{lem}
\begin{proof}
For $\alpha\geq\beta$, this is straightforward. Also, since
(\ref{our_eigenop}) does not 
have continuous spectrum, the eigenvalues are
continuous functions of the parameters. Hence, it is sufficient to show
that zero is not an eigenvalue of (\ref{our_eigenop}) when
$\alpha<\beta$.

Thus, letting $\lambda=0$, and $\xi\bydef\beta-\alpha>0$ in
(\ref{our_eigenop}) yields
\[
\varphi''-\left[\frac{1}{4}\left(\beta-\xi x\right)^2
 -\frac{\xi}{2}\right]\varphi=0,\quad 
 \varphi(0)=\varphi(1)=0.
\]
which can be further transformed by letting
$x=\xi^{-1}\beta+\sqrt{2}\xi^{-1/2}y$ in
\begin{equation}
\varphi''-(y^2-1)\varphi=0,\quad \varphi(A)=\varphi(A+B)=0,
\label{red_eigenop}
\end{equation}
where
\[
A=-\frac{\sqrt{2}}{2}\xi^{-1/2}\beta\quad\text{and}\quad B=\frac{\sqrt{2}}{2}\xi^{1/2}.
\]
The general solution to (\ref{red_eigenop}) is given by
\[
\varphi(y)=\e^{-y^2/2}\left(c_1+c_2\int_0^y\e^{s^2}\d s\right)
\]
On applying the boundary conditions, we see that a nontrivial solution
exists if, and only if, we have
\[
0=\int_0^{A+B}\e^{s^2}\d s - \int_0^{A}\e^{s^2}\d s =\int_A^{A+B}\e^{s^2}\d s.
\]
The last equality and the positiveness of the integrand implies
$B=0$, and hence $\xi=0$.
\end{proof}

\begin{prop}
The initial value problem defined by Equation~(\ref{sawrk_pde}) and
$w(0,x)=w_0(x)$, with $w_0\in L^1([0,1])$ is well posed and $w(t,x)\in
C^\infty(\mathbb{R}^+\times[0,1])$. Furthermore, we must have
\[
\lim_{t\to\infty}w(t,x)=0,\quad x\in[0,1].
\]
\end{prop}
\begin{proof}
Let $\varphi_j$ satisfy (\ref{our_eigenop}) with $\lambda_j$. Given
$f\in L^2([0,1])$ we set 
\[
f=\sum_{j\geq0}\hat{f}(j)\varphi_j
\]
Also, as in \cite{taylor1}, define for $s\in\mathbb{R}$
\[
\mathcal{D}_s=\left\{v\in L^1([0,1])\left| \sum_{j\geq0}|\hat{v}(j)|^2\lambda_j^s<\infty\right.\right\}
\]
Now, let
\begin{equation}
w(t,x)=\sum_{j\geq0}\hat{w}_0(j)e^{-t\lambda_j}\varphi_j(x).
\label{our_solution}
\end{equation}
For $t>0$, it is clear that $w$ 
satisfies
(\ref{our_equation}). If $w_0\in\mathcal{D}_s$ for $s>1/2$, then we
have a classical solution. In any case, however, notice that
(\ref{our_solution}) implies that $w(t,x)\in
C^\infty(\mathbb{R}^+\times[0,1])$, and that
\[
\lim_{t\to\infty}w(t,x)=0,\quad x\in[0,1].
\]
\end{proof}

Furthermore we have
\begin{lem}
Assume that $w_0\in L_2([0,1])$ and let $I(t)=\int_0^1w^2(t,x)\,\d
x$.
Then, we have
\[
I(t)\leq I(0)\e^{-2\lambda_0t}.
\]
\end{lem}
\begin{proof}
From the Fourier representation of $w(t,x)$, we have that
\[
I(t)=\sum_{j=0}^{\infty}\hat{w}^2_0(j)\e^{-2\lambda_jt}\leq 
\sum_{j=0}^{\infty}\hat{w}^2_0(j)\e^{-2\lambda_0t}=I(0)\e^{-2\lambda_0t}.
\]
\end{proof}

The solution given by (\ref{our_solution}), while well defined and
quite regular, has a major drawback: it does not satisfy, in general, the
required conservation laws, as it can be checked by starting with a
positive initial condition, and hence with positive mass. But the
decaying property of the (\ref{our_solution}) implies that the mass
will go to zero as time goes to infinity.

We shall give up as little regularity as possible, and look for a
solution in the class $C^\infty(\mathbb{R}^+\times(0,1))$. Thus, we
shall write
\begin{equation}
p(t,x)=q(t,x)+p_{\mathrm{D}}(t,x),\label{weak_form}
\end{equation}
where $q(t,x)$ satisfies (\ref{wrk_pde}) without boundary conditions,
and $p_{\mathrm{D}}(t,x)$ is a distribution solution with support in
$(0,\infty)\times\{0,1\}$. In 
this case, we must have, for some pair of nonnegative integers $M$ and
$M'$ that
\begin{equation}
p_{\mathrm{D}}(t,x)=\sum_{k=0}^Ma_k(t)\delta_0^k+ \sum_{k=0}^{M'}b_k(t)\delta_1^k,
\label{must_be}
\end{equation}
where $\delta_{x_0}^k$ means the $k$-th derivative of the delta
distribution at $x_0$.

Before proceeding, we must indicate precisely what we mean by a
weak solution in this case.

\begin{deff}\label{our_wsd}
A weak solution to (\ref{wrk_pde}) will be a distribution with support
in $[0,1]$ that satisfies
\begin{align*}
-\int_0^\infty\int_0^1p(t,x)\partial_t\phi(t,x)\d x\d t=&
\int_0^\infty\int_0^1p(t,x)\left[x(1-x)\partial^2_{x}\phi(t,x)+\right.\\
&\left.+x(1-x)(\beta+(\alpha-\beta)x)\partial_x\phi(t,x)\right]\d x\d t +\\
&+\int_0^1p^0(x)\phi(0,x)\d x,
\end{align*}
where 
\[
\phi(t,x)\in C^\infty_c\left([0,\infty)\times[0,1]\right).
\]
\end{deff}
\begin{rmk} Notice that the test functions in definition~\ref{our_wsd} are
  required to be of compact support in $[0,1]$ and not just in $(0,1)$
  as usual. Similar definitions have been given in other contexts; see
  for instance \cite{llx}. 

This definition can be recasted in the framework   of usual
  distribution theory, by introducing the compactly supported
  distribution  
\[
u(t,x)=
\sum_{k=0}^Ma_k(t)\delta_0^k+\sum_{k=0}^{M'}b_k(t)\delta_1^k +
\chi_{[0,1]}(x)q(t,x), 
\]
where $\chi_{[0,1]}$ is the characteristic function of unit interval. In
this case, the distribution can act in $C^\infty(\mathbb{R})$ and its
entirely determined by the behavior in the support; see for instance
\cite{hormander}. We shall abuse language and shall, henceforth,
identify $u(t,x)$ with $p(t,x)$.
\end{rmk}
We now can state the following important result:
\begin{lem}
\begin{enumerate}
\item Given $p^0(x)\in L^1([0,1])$, there is a unique weak solution
  $p(t,x)$ of (\ref{wrk_pde}) such that $p(t,x)\in
  C^\infty(\mathbb{R}^+\times(0,1))$ that satisfies
\[
\frac{\d}{\d t}\int_0^1p(t,x)\d x = 0\quad\text{and}\quad
\frac{\d}{\d t}\int_0^1\psi(x)p(t,x)\d x = 0.
\] 
\item This unique solution can be written as 
\[
p(t,x)=q(t,x)+a(t)\delta_0(x)+b(t)\delta_1(x),
\]
with
\[
a(t)=\int_0^t q(s,0)\d s  \quad\text{and}\quad
b(t)=\int_0^t  q(s,1)\d s ,
\]
where $q(t,x)$ is given by (\ref{our_solution}).
\end{enumerate}
\end{lem}
\begin{proof}
We begin by substituting (\ref{weak_form}), with $p_{\mathrm{D}}$ given by
(\ref{must_be}) into Definition~\ref{our_wsd} to obtain
\begin{align*}
&-(1)^{k+1}\int_0^\infty\left\{%
\sum_{k=0}^M a_k(t)\partial_t\partial_x^k\phi(t,0)+%
\sum_{k=0}^{M'} b_k(t)\partial_t\partial_x^k\phi(t,1)+\right\}\d t\\
&-\int_0^\infty\int_0^1q(t,x)\partial_t\phi(t,x)\d t=\\
&=\int_0^\infty\sum_{k=0}^M
a_k(t)(-1)^k\sum_{j=0}^{\max(k,2)}\binom{j}{k}\partial_x^j\left[x(1-x)\right]|_{x=0}\partial_x^{k-j+2}\phi(t,0)\d
t+ \\ 
&+\int_0^\infty\sum_{k=0}^{M'}
b_k(t)(-1)^k\sum_{j=0}^{\max(k,2)}\binom{j}{k}\partial_x^j\left[x(1-x)\right]|_{x=1}\partial_x^{k-j+2}\phi(t,1)+\d
t\\ 
&+\int_0^\infty\int_0^1q(t,x)x(1-x)\partial_x^2\phi(t,x)\d x\d t+\\
&+\int_0^\infty\sum_{k=0}^M
a_k(t)(-1)^k\sum_{j=0}^{\max(k,3)}\binom{j}{k}\partial_x^j\left[x(1-x)(\beta+\eta
  x)\right]|_{x=0}\partial_x^{k-j+1}\phi(t,0)\d
t+ \\ 
&+\int_0^\infty\sum_{k=0}^{M'}
b_k(t)(-1)^k\sum_{j=0}^{\max(k,3)}\binom{j}{k}\partial_x^j\left[x(1-x)(\beta+\eta
  x)\right]|_{x=1}\partial_x^{k-j+1}\phi(t,1)\d
t+ \\ 
&+\int_0^\infty\int_0^1q(t,x)x(1-x)(\beta+\eta x)\partial_x\phi(t,x)\d x\d t+\\
&+\int_0^1p^0(x)\phi(0,x)\d x.
\end{align*}
In the calculation above, we used that
\[
\partial^i_x\left[x(1-x)\right]=0,\quad i>2\quad\text{and}\quad
\partial^i_x\left[x(1-x)(\beta+\eta x)\right]=0,\quad i>3,
\]
Using that $q$ is smooth, integrating by parts, and using
(\ref{wrk_pde}  yields the following:
\begin{align*}
&(-1)^{k+1}\int_0^\infty\left\{%
\sum_{k=0}^M a_k(t)\partial_t\partial_x^k\phi(t,0)+%
\sum_{k=0}^{M'} b_k(t)\partial_t\partial_x^k\phi(t,1)+\right\}\d t=\\
&=\int_0^\infty\left[q(t,1)\phi(t,1)+q(t,0)\phi(t,0)\right]\d t+\\
&+\int_0^\infty\sum_{k=0}^M
a_k(t)(-1)^k\sum_{j=0}^{\max(k,2)}\binom{j}{k}\partial_x^j\left[x(1-x)\right]|_{x=0}\partial_x^{k-j+2}\phi(t,0)\d
t+ \\ 
&+\int_0^\infty\sum_{k=0}^{M'}
b_k(t)(-1)^k\sum_{j=0}^{\max(k,2)}\binom{j}{k}\partial_x^j\left[x(1-x)\right]|_{x=1}\partial_x^{k-j+2}\phi(t,1)+\d
t\\ 
&+\int_0^\infty\sum_{k=0}^M
a_k(t)(-1)^k\sum_{j=0}^{\max(k,3)}\binom{j}{k}\partial_x^j\left[x(1-x)(\beta+\eta
  x)\right]|_{x=0}\partial_x^{k-j+1}\phi(t,0)\d
t+ \\ 
&+\int_0^\infty\sum_{k=0}^{M'}
b_k(t)(-1)^k\sum_{j=0}^{\max(k,3)}\binom{j}{k}\partial_x^j\left[x(1-x)(\beta+\eta
  x)\right]|_{x=1}\partial_x^{k-j+1}\phi(t,1)\d
t.
\end{align*}
First, we look at $x=0$. Since the above must hold for any test
function we must have, for $k=0,1$, that
\begin{align*}
-\int_0^\infty a_0(t)\partial_t\phi(t,0)\d t&=
\int_0^\infty  q(t,0)\phi(t,0)\d t\\
\int_0^\infty a_1(t)\partial_t\partial_x\phi(t,0)\d t&=
\int_0^\infty\sum_{l=0}^3
a_l(t)(-1)^l\partial_x^l\left[x(1-x)(\beta+\eta
x)\right]|_{x=0}\partial_x\phi(t,0)\d t\\
\end{align*}
For $2\leq k\leq M$, we have
\begin{align*}
&(-1)^{k+1}\int_0^\infty a_k(t)\partial_t\partial_x^k\phi(t,0)\d t=\\
&=\int_0^\infty\sum_{l=k-2}^k
a_l(t)(-1)^l\binom{l-(k-2)}{l}\partial_x^{l-(k-2)}\left[x(1-x)\right]|_{x=0}\partial_x^k\phi(t,0)\d
t+\\
&+\int_0^\infty\sum_{l=k-1}^{\min(k+2,M)}
a_l(t)(-1)^l\binom{l-(k-1)}{l}\partial_x^{l-(k-1)}\left[x(1-x)(\beta+\eta
  x)\right]|_{x=0}\partial_x^k\phi(t,0)\d t.\\
\end{align*}
For $k=M+1,M+2$, we find:
\begin{align*}
&0=\\
&=\int_0^\infty\sum_{l=k-2}^M
a_l(t)(-1)^l\binom{l-(k-2)}{l}\partial_x^{l-(k-2)}\left[x(1-x)\right]|_{x=0}\partial_x^k\phi(t,0)\d
t+\\
&+\int_0^\infty\sum_{l=k-1}^M
a_l(t)(-1)^l\binom{l-(k-1)}{l}\partial_x^{l-(k-1)}\left[x(1-x)(\beta+\eta
  x)\right]|_{x=0}\partial_x^k\phi(t,0)\d t.\\
\end{align*}
For $k=M+2$, the relation above is identically zero but,
for $k=M+1$, we have that
\[
0=(-1)^MM\int_0^\infty a_M(t)\partial_x^{M+1}\phi(t,0)\d t
\]
Hence $a_M(t)\equiv0$.

Considering $k=M$, yields
\begin{align*}
&(-1)^{M+1}\int_0^\infty a_{M}(t)\partial_t\partial^{M}_x\phi(t,0)\d t=\\
&\int_0^\infty\left((-1)^{M-1}(M-1)
a_{M-1}(t)+(-1)^MM(M-1)a_M(t)+(-1)^MM\beta
a_M(t)\right)\partial_x^M\phi(t,0)\d t. 
\end{align*}
Since $a_M(t)\equiv0$, we have that $a_{M-1}(t)\equiv0$ as well.
For $k=M-1$, we have that
\begin{align*}
&(-1)^{M}\int_0^\infty
  a_{M-1}(t)\partial_t\partial^{M-1}_x\phi(t,0)\d t=\\
&\int_0^\infty\left((-1)^{M-2}(M-2)a_{M-2}(t) 
-(-1)^{M-1}(M-1)(M-2)a_{M-1}(t)\right)\partial_x^M\phi(t,0)\d
  t\\+
&\int_0^\infty\left((-1)^{M-1}(M-1)\beta
  a_{M-1}(t)+(-1)^MM(M-1)(\eta-\beta)a_{M}(t)\right)\partial_x^M\phi(t,0)\d t.
\end{align*}
Again,we have $a_M(t)\equiv a_{M-1}(t)\equiv0$; thus
$a_{M-2}(t)\equiv0$.

For $1\leq k\leq M-2$, we have a linear relation involving $a_i(t)$,
$i=k,\ldots,k+3$ (when $k=1$, we have $i=1,\ldots,3$). If three of
them are zero, then the remaining one is also zero. Thus, starting
with $k=M-2$ and proceeding inductively, we find that $a_k(t)\equiv0$
for $k=1,\ldots,M$. Therefore, only $a_0(t)$ can be  nonzero. 

An analogous argument shows also that only $b_0(t)$ can be nonzero as
well. We now drop the subscripts and determine their values.

Integrating by parts, the corresponding relation for $a(t)$, we obtain 
\[
\int_0^\infty a(t)\partial_t\phi(t,0)=
-\int_0^\infty  \int_0^tq(s,0)\d s\phi(t,0)
\]
Hence
\[
a(t)=\int_0^tq(s,0)\d s +a_0.
\]
A similar calculation shows that
\[
b(t)=\int_0^tq(s,1)\d s +b_0,
\]
It remains only to show that the conservation laws are
satisfied. Substituting the found solution on them, we find
\[
a'(t)+b'(t)-q(t,1)-q(t,0)=0\quad\text{and}\quad
a'(t)-q(t,1)=0
\]
respectively, which are obviously satisfied.
\end{proof}

\bibliography{moran}
\bibliographystyle{abbrv}

\end{document}